\documentclass{article}
\usepackage{graphicx} % Required for inserting images

\usepackage{comment}
\usepackage{float}
\usepackage{amsmath}
\usepackage{listings}
\usepackage{amsfonts}
\usepackage{pgf,tikz,pgfplots}
\pgfplotsset{compat=1.15}
\usepackage{mathrsfs}
\usepackage{comment}
\usetikzlibrary{arrows}
\usepackage[normalem]{ulem}
\usepackage{amsmath,amssymb,amsthm}
\usepackage{subcaption} 
\usepackage{booktabs}

\usepackage{cite}
\usepackage{authblk}
\usepackage[margin=3.5cm]{geometry} % so tables are fully included in the text zone

\def\Xint#1{\mathchoice
{\XXint\displaystyle\textstyle{#1}}%
{\XXint\textstyle\scriptstyle{#1}}%
{\XXint\scriptstyle\scriptscriptstyle{#1}}%
{\XXint\scriptscriptstyle\scriptscriptstyle{#1}}%
\!\int}
\def\XXint#1#2#3{{\setbox0=\hbox{$#1{#2#3}{\int}$}
\vcenter{\hbox{$#2#3$}}\kern-.5\wd0}}

\def\dashint{\Xint-}

\newtheorem{theorem}{Theorem}
\newtheorem{proposition}{Proposition}
{\theoremstyle{remark} \newtheorem{remark}{Remark}}

\title{A Godunov-type in-cell discontinuous reconstruction shock-capturing method for the Self-Organized Hydrodynamics model}

\author[1]{Christophe Berthon}
\author[1]{Marie Compain}
\author[2]{Carlos Par\'es}
\author[3]{Ernesto Pimentel-Garc{\'\i}a}

\affil[1]{Laboratoire de Math\'ematiques Jean Leray, CNRS UMR 6629, Nantes Universit\'e, 2 rue de la Houssini\`re, BP 92208, 44322 Nantes, France}
\affil[2]{Departamento de An\'alisis Matem\'atico, Estad{\'\i}stica e Investigaci\'on Operativa, y Matem\'atica Aplicada, Universidad de M\'alaga, Campus de Teatinos s/n, M\'alaga 29080, Spain}
\affil[3]{Departamento de Matem\'atica Aplicada, Universidad de M\'alaga, Campus de Teatinos s/n, M\'alaga 29080, Spain}

\date{\today}

\begin{document}

\maketitle

\begin{abstract}
The present paper is devoted to show the correct convergence of a numerical method based on Godunov and the in-cell discontinuous reconstruction technique for the Self-Organized Hydrodynamics model. The $2\times2$
system under consideration is strictly hyperbolic but never recasts in
a conservation form. With ambiguous nonconservative products, we adopt
a non usual Generalized Rankine-Hugoniot jump condition to suitably
define the shock wave solutions. In a previous work the authors determine the
simple waves of the nonconservative system studied here and they are
able to construct the solutions of the Riemann problem of interest. These solutions are used to obtain a Godunov-type scheme and it is combined with an in-cell discontinuous reconstruction technique to obtain the correct convergence to the exact solution even in the presence of shocks.
\end{abstract}

%%%%%%%%%%%%%%%%%%%%%%%%%%%%%%%%%%%%%%%%%%%%%%%%%%%%%%%%%%%
%%%%%%%%%%%%%%%%%%%%%%%%%%%%%%%%%%%%%%%%%%%%%%%%%%%%%%%%%%%
%%%%%%%%%%%%%%%%%%%%%%%%%%%%%%%%%%%%%%%%%%%%%%%%%%%%%%%%%%%

\section{Introduction}
In this work, we develop a numerical scheme capable of converging to the physical solutions of a nonconservative hyperbolic PDE system describing self-organized hydrodynamics. The system was derived in \cite{Degond2012,Degond2008} as the hydrodynamic limit of the classical Cucker--Smale model \cite{vicsek1995novel}, and describes the collective motion of self-propelled particles through alignment interactions.

Such models have been applied to a wide range of collective phenomena. In biology, they describe collective cell migration, tissue formation, and biochemical signal propagation \cite{Creppy2016}. In crowd dynamics, they model pedestrian flows, congestion, and evacuation processes \cite{Degond2013}. They have also been used to study opinion dynamics and social self-organisation \cite{Degond2017a}, as well as self-organising systems in robotics, nanotechnology, and autonomous drone fleets \cite{Degond2014,Degond2025,Degond2017b}. In ecology, they provide models for collective behaviour in aquatic ecosystems, such as fish schools and plankton colonies \cite{Degond2017}.

The one-dimensional Self-Organized Hydrodynamics (SOH) model considered in this work was proposed in \cite{Degond2008} (see also \cite{Degond2012,motsch2011numerical,Frouvelle2012}). It describes the evolution of the density $\rho$ and the velocity $(u,v)$ of self-propelled particles. Since all particles have the same constant speed, we impose
\begin{equation}
\sqrt{u^2+v^2}=1.
\end{equation}
Introducing the velocity direction through $u=\cos(\theta)$ and $v=\sin(\theta)$, the model becomes
\begin{equation}\label{pde}
\left\{
\begin{array}{l}
\displaystyle \partial_t \rho + \partial_x(\rho\cos(\theta)) = 0, \\
\displaystyle \partial_t\theta + \partial_x(c\sin(\theta)) -\lambda\frac{\sin(\theta)}{\rho}\partial_x\rho = 0.
\end{array}
\right.
\end{equation}
where $c\in (0,1)$ and $\lambda > 0$ are two constants prescribed by the physics (see \cite{motsch2011numerical,Degond2008}).
This nonlinear hyperbolic system is non-conservative, and therefore discontinuous solutions require a suitable definition of the non-conservative products across shock waves.

To derive the corresponding jump conditions, we consider the following viscous regularization introduced in \cite{Degond2012}
\begin{equation}\label{sys_visquous}
\left\{
\begin{array}{l}
\displaystyle \partial_t\rho + \partial_x(\rho\cos(\theta)) = 0, \\
\displaystyle \partial_t\theta + \partial_x(c\sin(\theta)) -\lambda\frac{\sin(\theta)}{\rho}\partial_x\rho = c_v\left(\frac{2}{\rho}\partial_x\theta\partial_x\rho +\partial_{xx}\theta\right),
\end{array}
\right.
\end{equation}
with $c_v>0$ a constant.
Following \cite{Sainsaulieu1991,compain:hal-05490904}, travelling-wave solutions of this viscous system lead to the generalized Rankine--Hugoniot conditions associated with \eqref{pde}: a pair of states $U_L = [\rho_L, \theta_L]^T$, $U_R = [\rho_R, \theta_R]^T$ defines a shock wave of speed $\sigma$ compatible with the viscous model \eqref{sys_visquous} as soon as the following conditions hold:
\begin{itemize}
\item[(SW1)] The states $U_L$, $U_R$, and the velocity $\sigma$ satisfy
\begin{equation}\label{RH}
\left\{
\begin{array}{l}
\displaystyle \sigma(\rho_R - \rho_L) = \rho_R \cos(\theta_R) - \rho_L \cos(\theta_L), \\
\displaystyle G_\sigma (\theta_R) - G_\sigma (\theta_L) = 0,
\end{array}
\right.
\end{equation}
where $G_\sigma$ is given for $\cos(\theta) \neq \sigma$ by
\begin{equation}\label{G_sigma}
G_\sigma'(\theta) = \frac{- \sigma + c \cos(\theta) - \frac{\lambda \sin(\theta)^2}{\cos(\theta) - \sigma}}{(\cos(\theta) - \sigma)^2}.
\end{equation}

\item[(SW2)] If $\theta_L < \theta_R$, for all $\theta \in (\theta_L, \theta_R)$, we have $G_\sigma(\theta) - G_\sigma(\theta_L) > 0$; \\
if $\theta_L > \theta_R$, for all $\theta \in (\theta_R, \theta_L)$, we have $G_\sigma(\theta) - G_\sigma(\theta_L) < 0$.

\item[(SW3)] For all $\theta$ between $\theta_L$ and $\theta_R$, we have $\cos(\theta) \neq \sigma$.
\end{itemize}

The function $G_\sigma(\theta)$ can be computed with elementary functions and its expression is detailed in Appendix \ref{appendix_1}. It is worth noting that these nonstandard jump conditions provide a remarkable counterexample to the equivalence between the evanescent viscosity approach \cite{Sainsaulieu1991} and the path theory \cite{DMLFM1995,lefloch1989,berthon2012}. Let us consider a nonconservative system under the following quasi-linear form:
\begin{equation}\label{sbl}
U_t + F(U)_x + \mathcal{B}(U)U_x = 0,
\end{equation}
which can also be reformulated as
\begin{equation} \label{eq:quasiA}
    U_t + \mathcal{A}(U)U_x = 0.
\end{equation}
In the path theory framework, non-conservative products are defined via a family of paths $\Phi(s; U_L, U_R)$, leading to jump conditions of the form
\begin{eqnarray}\label{RH2}
& & \sigma(U_R - U_L)  =  \int_0^1 \mathcal{A}(\Phi(s; U_L, U_R))\partial_s \Phi(s; U_L, U_R) \, ds\\
& & \qquad = F(U_R) - F(U_L) +  \int_0 ^1  \mathcal{B}(\Phi(s; U_L, U_R))\partial_s \Phi(s; U_L, U_R) \, ds.\nonumber
\end{eqnarray}

As demonstrated in \cite{Vicsek_Riemann}, the following negative result holds:

\begin{proposition} \label{prop:nopath}
It is not possible to find a family of paths that depends continuously on the parameter $\lambda \in [0, \infty)$ such that the jump condition \eqref{RH2} is equivalent to \eqref{RH} for all $\lambda$.
\end{proposition}

Observe that  \eqref{RH} can be written in a similar form to \eqref{RH2} as follows
\begin{equation}\label{RHtheta}
\sigma (U_R - U_L) = F(U_R) - F(U_L) + \int_{\theta_L}^{\theta_R} \mathcal{G}_\sigma(\theta) \, d\theta,
\end{equation}
where
$$
\mathcal{G}_\sigma(\theta) = 
\left[
\begin{array}{c}
0\\
\displaystyle
\frac{(\sigma - c \cos(\theta))(1  -(\cos(\theta) - \sigma)^2) +  \frac{\lambda \sin(\theta)^2}{\cos(\theta) - \sigma}}{(\cos(\theta) - \sigma)^2}
\end{array}
\right].$$
The integral in the right-hand side can be understood then as the contribution of the non-conservative product $\mathcal{B}(U)U_x$ at an admissible shock wave. 

Given an admissible weak solution $U(x,t)$ of \eqref{pde}, with bounded total variation, let us define the generalized integral
\begin{eqnarray}\label{dashint}
\dashint_{a}^{b} \mathcal{A}(U(x,t))U_x(x,t) dx & = &  \int_{a}^{b} \mathcal{A}(U(x,t))U_x(x,t) dx \nonumber \\
&+ &\sum_{k} \left(F(U_k^+) - F(U_k^-) +
\int_{\theta_k^-}^{\theta_k^+} \mathcal{G}_\sigma(\theta) \, d\theta \right),
\end{eqnarray}
where  the derivative in the integral at the right-hand side represents the standard one, that exists almost everywhere, and $(U_k^-, U_k^+)$ in the sum represent the left and right states at the $k$th discontinuity in the interval $[a,b]$.
Then, the following equality holds true
\begin{equation}\label{Riemannprop}
\int_a^b U(x,t_1) \, dx = 
 \int_a^b U(x,t_0) \, dx  +  \int_{t_0}^{t_1}\dashint_{a}^{b} \mathcal{A}(U(x,t))U_x(x,t) \,dx \,dt,
\end{equation}
for every space and time intervals $[a,b]$ and $[t_0, t_1]$.
\begin{remark} More formally, the weak solutions of the system according to the jump condition \eqref{RHtheta} satisfy the PDE \eqref{sbl} in the sense of Borel measures. The dashed integral \eqref{dashint} represents the measure of the interval $[a,b]$ at time $t$ according to the definition of $A(U)U_x$ as a Borel measure. 
\end{remark}

The objective of this work is to apply the methodology developed in \cite{chalons2020path, pimentelgarcia2022cell, pimentel2024drmood} to converge to the correct numerical solutions even in the presence of discontinuities but without the path-conservative framework. The idea is to combine a Godunov scheme with an in-cell discontinuous reconstruction based on the information given by the solutions of the Riemann problems that were studied in \cite{Vicsek_Riemann}.

To address such an objective, the
paper is organized as follows. In the next section, we recall the
algebra of the system and we show the ingredients to obtain the solution of the Riemann problems given in \cite{Vicsek_Riemann}. Next, in Section 3, we detail the Godunov scheme based on the Riemann solutions found in \cite{Vicsek_Riemann}. As in the case of systems of conservation laws, the Godunov method is obtained by averaging, at time $t^{n+1}$, the exact solutions to the Riemann problems associated with each intercell at time $t^n$. The averaging process introduces numerical viscosity into the method. Due to this, for non-conservative systems, the numerical approximations may not converge to the correct weak solutions of the problem: see \cite{castro2008}. In order to avoid the addition of numerical viscosity close to a shock wave, the in-cell reconstruction technique is applied in Section \ref{sec:in_cell_method}: first, a semi-discrete form of the scheme is presented, then the in-cell discontinuous reconstruction is explained in different subsections and combined with the Godunov scheme. Finally, in Section \ref{sec:numerical_tests} some numerical methods are shown to see the convergence of the new method to the exact solutions
and the failure of convergence of Godunov method.

%%%%%%%%%%%%%%%%%%%%%%%%%%%%%%%%%%%%%%%%%%%%%%%%%%%%%%%%%%%
%%%%%%%%%%%%%%%%%%%%%%%%%%%%%%%%%%%%%%%%%%%%%%%%%%%%%%%%%%%
%%%%%%%%%%%%%%%%%%%%%%%%%%%%%%%%%%%%%%%%%%%%%%%%%%%%%%%%%%%

\section{Main nature of the SOH model}
For the sake of completeness, we now present the main results
concerning the algebra of the evolution law \eqref{pde}. The
computations presented in this section were performed in
\cite{motsch2011numerical} and they are essential in the sequel of
this paper. In order to clarify the notations, we set
\begin{equation}\label{U}
U = \left[ \begin{array}{c} \rho \\ \theta \end{array} \right] \in [0,\infty)\times(-\pi, \pi],
\end{equation}
to write the system \eqref{pde} in the quasi-linear system (\ref{sbl}) with
$$
F(U) = \left[ \begin{array}{c} \rho \cos(\theta) \\ c \sin(\theta) \end{array} \right], 
\quad \mathcal{B}(U) = \left[ \begin{array}{cc} 0 & 0  \\\displaystyle -\lambda\frac{\sin(\theta)}{\rho} & 0 \end{array} \right].$$
It can also be written in quasilinear form {(\ref{eq:quasiA}) where we have introduced
\begin{equation}\label{A}
\mathcal{A}(U) = \left[ \begin{array}{cc} \cos \theta & -\rho \sin \theta \\ \displaystyle- \lambda \frac{\sin \theta}{\rho} & c \cos \theta
\end{array} \right].
\end{equation}

The above matrix admits two real eigenvalues given by
\begin{equation}\label{av1}
\mu_\pm(\theta) = \frac{c + 1}{2} \cos \theta \pm \gamma, \quad \gamma = \sqrt{\frac{(c-1)^2}{4} \cos^2 \theta + \lambda \sin^2 \theta}.
\end{equation}
Since $c<1$, the eigenvalues are distinct and then $\mathcal{A}(U)$ is
diagonalizable in ${\mathbb R}$. As a consequence, the system
\eqref{pde} is strictly hyperbolic $(0, \infty) \times (-\pi, \pi]$. Moreover, let us note that the
associated eigenvectors read
$$
R_- (\theta) = \left[ \begin{array}{c} \displaystyle c \cos \theta - \mu_-(\theta)\\  \frac{\lambda}{\rho} \sin
    \theta\end{array}
  \right] 
\quad \text{ and }
R_+ (\theta) = \left[ \begin{array}{c} \rho \sin
    \theta\\ \displaystyle \cos \theta - \mu_+(\theta) \end{array}
  \right].
$$

The characteristic field $(\mu_+,R_+)$ is neither genuinely nonlinear
nor linearly degenerated since $\nabla_U\mu_+\cdot R_+=0$ as soon as
$\theta=\pm\theta^\star$ with
\begin{equation}\label{eq:theta_star}
\theta^* = \arctan\left(\sqrt{\frac{4\lambda^2 - (2\lambda + c)(c-1)^2 }{\lambda (c +1)^2}}\right).
\end{equation}
Let us remark that the physical values of $c$ and $\lambda$ enforce
$4\lambda^2 - (2\lambda + c)(c-1)^2>0$ so that
$\theta^\star\in(-\pi,\pi]$ exists. In addition, we notice that
  $\mu_+(\theta)$ reaches its maximum for $\pm\theta^\star$ and its
  minimum for $\theta=0,\pi$.

Likewise, the characteristic field $(\mu_-,R_-)$ is neither genuinely
nonlinear nor linearly degenerated because $\nabla_U\mu_-\cdot R_-$
vanishes for $\theta=\pm(\theta^\star-\pi)$. Moreover, we notice that
the maximum value of $\mu-(\theta)$ is reached for
$\theta=\pm(\theta^\star-\pi)$ while its minimum values are obtained
for $\theta=0,\pi$.

For the sake of simplicity in the notations, we set
$$
\overline{\mu}_-= \mu_-(\pi-\theta^*) = \min(\mu_-(\theta)), \quad \overline{\mu}_+= \mu_+(\theta^*) = \max(\mu_+(\theta)).
$$

To illustrate these behaviors, in Figure \ref{fig:eigenvalues} we display the graphs of the
eigenvalues and the values of $\theta$ corresponding to the change of
nature of the characteristic fields; formally the values of $\theta$
where the genuinely nonlinear character is lost.

\begin{figure}[htpb]
\begin{center}
\includegraphics[width=0.6\linewidth]{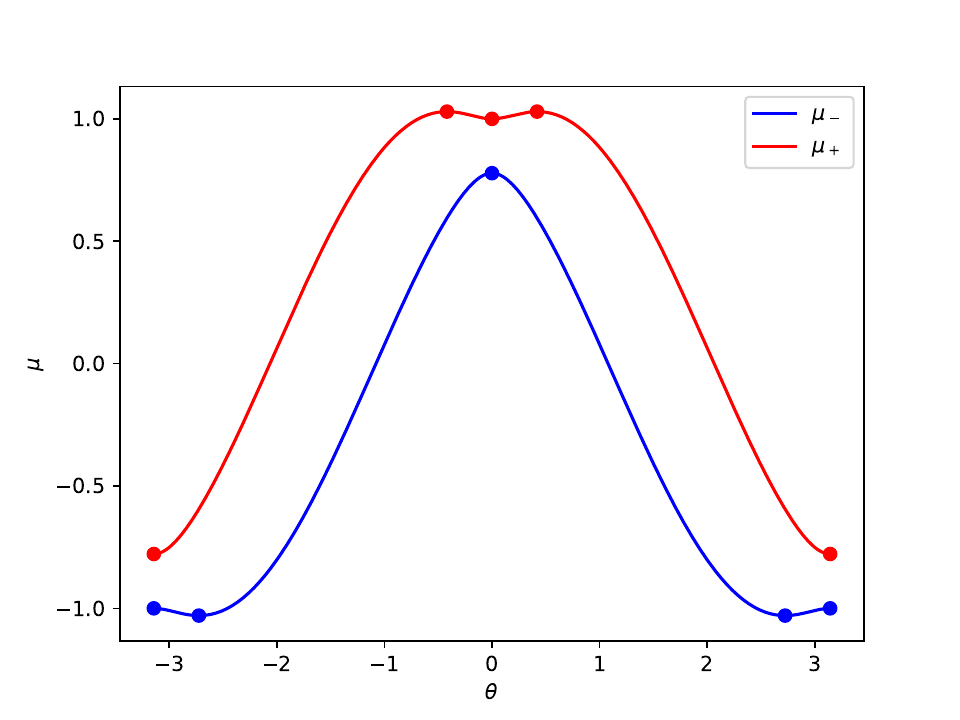}
\caption{Graphs of the eigenvalues $\mu_\pm(\theta)$ for $c = 0.8$ and $\lambda = 0.2$. The angles $\pm\pi$, $0$, $\theta^* - \pi$, $ \pi-\theta^*$ (resp. $\pm\pi$, $0$, $\pm\theta^*$) in which the first (resp. second) characteristic field is not genuinely nonlinear,  are marked with a green circle (resp. a red star) in the corresponding graph.}
\label{fig:eigenvalues}
\end{center}
\end{figure}

Because of the complexity of the involved algebra calculus, in the
present work we decide to avoid eventual change in the nature of the
characteristic fields. As a consequence, we focus here in the analysis
and numerical resolution of system \eqref{pde} in the region
\begin{equation}\label{defOmega} 
\Omega = (0, \infty) \times (\theta^*, \pi - \theta^*),
\end{equation}
where
$\nabla_U\mu_\pm\cdot R_\pm\not=0$. With some abuse in the definition,
we then enforce our study in regions where the characteristic fields
are genuinely nonlinear. However, it can be easily checked that, if $U(x,t) =
[\rho(x,t), \theta(x,t)]^T$ is a solution of \eqref{pde} such that
$U(x,t) \in \Omega$ for all $x$, $t$, then $\tilde U(x,t) =
[\rho(-x,t), -\theta(-x,t)]^T$ is a solution satisfying $\widetilde
U(x,t) \in \widetilde \Omega: = (0, \infty) \times (\theta^*-\pi,
-\theta^*) $ so that all the conclusions of our analysis also apply to
$\widetilde \Omega$.

Since it will be useful in the sequel, this section is concluded by
the definition of the inverse function of $\mu_\pm$. We then define the functions
\begin{equation}\label{invav}
\chi_\pm(x) = \arccos\left(\frac{x(c+1) \pm \sqrt{x^2(c-1)^2 + 4\lambda(c + \lambda - x^2)}}{2(c+ \lambda)}\right).
\end{equation}
Hence, it can be checked that
\begin{itemize}
 \item $\chi_+:[\mu_-(\pi-\theta^*), \mu_-(0)] \to \mathbb{R}$ is the inverse function of $\mu_-$ in the interval $[0, \pi - \theta^*]$;
\item $\chi_-:[\mu_+(\pi), \mu_+(\theta^*)] \to \mathbb{R}$ is the inverse function of $\mu_+$ in the interval $[\theta^*, \pi]$.
\end{itemize}

\section{Riemann problems and Godunov method} \label{sec:Godunov}

This section is devoted to the development of a Godunov method for the SOH model (\ref{pde}).
In \cite{Vicsek_Riemann},  the solution of the Riemann problem associated to \eqref{pde} with initial conditions
\begin{equation}\label{eq:Riemann_general}
U(x,0 )  = \begin{cases}
    U_L & \text{if $x < 0$;} \\
    U_R & \text{ if $x > 0$;}
\end{cases}
\end{equation}
has been studied. It is shown that, for close enough states $U_L$, $U_R$, the Riemann problem has a unique self-similar solution
$$
U(x,t ) = \mathcal{R} \left( \frac{x}{t}; U_L, U_R\right),$$
consisting of two simple waves (shocks and/or rarefactions) that link the initial states with an intermediate state $U_m = [\rho_m, \theta_m]^T$. The structure of the solution depends on $\theta_m$:

\begin{enumerate}

\item If $\theta_L < \theta_m < \theta_R$ the solution of the Riemann problem consists of two shocks. 

\item If $\theta_R < \theta_m < \theta_L$ the solution of the Riemann problem consists of two rarefaction waves. 

\item If $\theta_m < \theta_L$ and $ \theta_m < \theta_R$ the solution of the Riemann problem consists of a 1-rarefaction and a 2-shock wave. 

\item If $\theta_m > \theta_L$ and $ \theta_m > \theta_R$ the solution of the Riemann problem consists of a 1-shock and a 2-rarefaction waves. 

\end{enumerate}

For completeness, the equation to be solved to find $\theta_m$ is detailed in Appendix \ref{appendix_2}. In what follows the notations
$$
U_0 = U_L, \quad U_1 = U_m, \quad U_2 = U_R,$$
will be used to represent the three states involved by the solution of a Riemann problem.

Once the solution of the Riemann problem is known, the method of Godunov can be applied to the system in the usual form by averaging at time $t_{n+1}$ the solutions of the Riemann problems at the intercells at time $t_n$:
\begin{equation}\label{Godunov}
U_i^{n+1} = 
\frac{1}{\Delta x}\left(  
\int_{0}^{\Delta x/2} \mathcal{R}\left(\frac{x}{\Delta t}; U_{i-1}, U_i\right) \,dx
+ \int_{-\Delta x/2}^0 \mathcal{R}\left(\frac{x}{\Delta t}; U_{i}, U_{i+1}\right) \,dx \right),
\end{equation}
where
\begin{itemize}
    \item $\Delta x$ and $\Delta t$ are the space and the time steps respectively. They are assumed to be constant in the description of the numerical method for simplicity.
    \item $I_{i} = \left[x_{i-\frac{1}{2}}, x_{i-\frac{1}{2}}\right]$ are the computational cells with length $\Delta x$.
    \item $t_n = n\Delta t$, $n=0,1,2,...$.
	\item $U_{i}^n \approx \displaystyle \frac{1}{\Delta x}\int_{x_{i+\frac{1}{2}}}^{x_{i-\frac{1}{2}}} U(x,t_n) \,dx$ is the approximation to the cell average of the solution at the $i$-th cell in time $t_n$.

\end{itemize}
The usual CFL-1/2 condition
\begin{equation}\label{CFL}
 \max_i\{ |\mu_\pm(U_i^n) |\}\frac{\Delta t}{\Delta x} \leq \frac{1}{2},
\end{equation}
is assumed so that the waves issue of the intercells at time $t_n$ do not reach the cell-center $x_i$ before time $t_{n+1}$.

Let us write method \eqref{Godunov} in an equivalent form that will be useful to introduce its in-cell discontinuous reconstruction extension. First, it can be easily checked that  \eqref{Godunov} can be rewritten as follows:
\begin{equation}\label{Godunov_1}
U_i^{n+1} = 
 U_i^n - \frac{\Delta t}{\Delta x} (D_{i-1/2}^+ + D_{i + 1/2}^-),
\end{equation}
where
\begin{eqnarray}
D_{i-1/2}^+ & = & -\frac{1}{\Delta t}\int_{0}^{\Delta x/2} \left(\mathcal{R}\left(\frac{x}{\Delta t}; U^n_{i-1}, U^n_i\right) - U^n_i\right)  \,dx, \label{D-}\\
% & = & -\int_{0}^{{\Delta x}/({2\Delta t})} \left(\mathcal{R}(s; U_{i-1}, U_i) - U_i\right)  \,ds; \nonumber \\
D_{i+1/2}^- & = & -\frac{1}{\Delta t}\int_{-\Delta x/2}^{0} \left(\mathcal{R}\left(\frac{x}{\Delta t}; U_{i}^n, U_{i+1}^n\right) - U^n_i\right)  \,dx.\label{D+} \\
% & = & -\int_{-\Delta x/({2\Delta t})}^0 \left(\mathcal{R}(s; U_{i1}, U_{i+1}) - U_i\right)  \,ds. \nonumber
\end{eqnarray}
Applying now \eqref{Riemannprop} in the space intervals $[0, \Delta x/2]$,
$[-\Delta x/2, 0]$ and the time interval $[0, \Delta t]$ we obtain
\begin{equation}
D_{i-1/2}^+ = D^+(U^n_{i-1}, U^n_i), \quad D_{i+1/2}^- = D^-(U^n_{i}, U^n_{i+1}),
\end{equation}
where the functions $D^\pm: \Omega \times \Omega \to \mathbb{R}^2$ are defined as follows
\begin{equation}
D^+(U_L, U_R)  =  -\frac{1}{\Delta t}
\int_0^{\Delta t}\dashint_{0}^{\Delta x/2} \mathcal{A}\left(\mathcal{R}\left(\frac{x}{t}; U_{L}, U_R\right)\right) \partial_x \mathcal{R}\left(\frac{x}{t}; U_L, U_R\right) \,dx dt,\label{D-2}
\end{equation}
\begin{equation}
D^-(U_L, U_R) =  -\frac{1}{\Delta t}
\int_0^{\Delta t}\dashint_{-\Delta x/2}^{0} \mathcal{A}\left(\mathcal{R}\left(
\frac{x}{t} ; U_{L}, U_{R} \right)\right) \partial_x \mathcal{R}\left(\frac{x}{t}; U_{L}, U_{R} \right)   \,dx dt,\label{D+2} 
\end{equation}
for all $U_L, U_R \in \Omega$.

Let us compute the fluctuations $D^\pm$ in some particular cases:
\begin{itemize}
\item If $U_L$, $U_R$ can be joined by a $j$-rarefaction wave $j = 1,2$, and $\mu(U_R) \leq  0$, where $\mu = \mu_-$ or $\mu = \mu_+$ is the corresponding eigenvalue, $D^+(U_L, U_R) = 0$; otherwise, for every $t$
$$
x \in [0, \Delta x/2] \to \mathcal{R}\left(\frac{x}{t}; U_{L}, U_R\right) \in \Omega,$$
is a parameterization of an arc of the integral curve of the corresponding characteristic field that contains both $U_L$ and $U_{R}$. If the eigenvalue $\mu$ is used to parameterize this arc
$$
\mu \in [\mu_0, \mu(U_R)] \to \widetilde{U}(\mu) \in \Omega,
$$
where
$$
\mu_0 = \mu(\mathcal{R}\left(0; U_{L}, U_R\right)),$$
the integral does not change and we obtain
%$$ -> compilation problem with this
\begin{align}
& \dashint_{0}^{\Delta x/2} \mathcal{A}\left(\mathcal{R}\left(\frac{x}{t}; U_{L}, U_R\right)\right) \partial_x \mathcal{R}\left(\frac{x}{t}; U_{L}, U_R\right) \,dx  \\
&= \int_{0}^{\Delta x/2} \mathcal{A}\left(\mathcal{R}\left(\frac{x}{t}; U_{L}, U_R\right)\right) \partial_x \mathcal{R}\left(\frac{x}{t}; U_{L}, U_R\right) \,dx \\
&= \int_{\mu_0}^{\mu(U_R)} \mathcal{A}(\widetilde{U}(\mu) )\widetilde{U}'(\mu) \, d\mu \\
&= \int_{\mu_0}^{\mu(U_R)} \mu \widetilde{U}'(\mu) \, d\mu \\
& = \mu(U_R)U_R - \mu_0 \mathcal{R}\left(0; U_{L}, U_R\right) - \int_{\mu_0}^{\mu(U_R)}  \widetilde{U}(\mu) \, d\mu ,
\end{align}
%$$
where an integration by parts has been performed in the last step.
Since the result is independent of $t$ we obtain
$$
D^+(U_L, U_R) =  \mu(U_R)U_R - \mu_0 \mathcal{R}\left(0; U_{L}, U_R\right) - \int_{\mu_0}^{\mu(U_R)}  \widetilde{U}(\mu) \, d\mu . $$

\item Similarly,  if $U_{L}$, $U_{R}$ can be joined by a $j$-rarefaction wave $j = 1,2$, and $\mu(U_L) \geq  0$, where $\mu$ is the corresponding eigenvalue, 
then $D^-(U_L, U_R) = 0$; otherwise
$$
D^-(U_L, U_R) =  \mu_0 \mathcal{R}\left(0; U_{L}, U_{R}\right) - \mu(U_L)U_L   - \int_{\mu(U_L)}^{\mu_0} U(\mu) \, d\mu . $$

\item If $U_{L}$, $U_R$ can be joined by a $j$-shock wave with velocity $\sigma$ and $\sigma \leq 0$ then $D^+(U_L, U_R) = 0$; otherwise we have
%$$
\begin{align*}
& \dashint_{0}^{\Delta x/2} \mathcal{A}\left(\mathcal{R}\left(\frac{x}{t}; U_{L}, U_R\right)\right) \partial_x \mathcal{R}\left(\frac{x}{t}; U_{L}, U_R\right) \,dx  \\
&= F(U_{R}) - F(U_{L}) +
\int_{\theta_L}^{\theta_{R}} \mathcal{G}_\sigma(\theta) \, d\theta \\
&= \sigma (U_R - U_L ),
\end{align*}
%$$
where \eqref{RHtheta} has been used. Since the result is independent of $t$ we obtain
$$
D^+(U_L, U_R) = \sigma (U_R - U_L ). $$

\item Similarly, if $U_{L}$, $U_{R}$ can be joined by a $j$-shock wave with velocity $\sigma$ and $\sigma \geq 0$ then $D^-(U_L, U_R) = 0$; otherwise 
$$D^-(U_L, U_R) = \sigma (U_R - U_L ). $$

    \end{itemize}

In the general case one has
$$
D^\pm(U_L, U_R)  = D_1^\pm(U_L, U_R) + D_2^\pm(U_L, U_R),$$
where $D_j^\pm(U_L, U_R)$, $j=1,2$ are the contributions of the $j$-waves, that are computed as follows:
\bigskip

\noindent \textbf{1-Wave Contributions ($D_1^-, D_1^+$)}:
\begin{itemize}
    \item \textbf{1-Shock ($\theta_L < \theta_m$)} with speed $\sigma_- = \displaystyle \frac{\rho_m \cos\theta_m - \rho_L \cos\theta_L}{\rho_m - \rho_L}$:
  \begin{itemize}
        \item If  $\displaystyle \sigma_- \le 0$:  
         $\displaystyle D_1^- = \sigma_- (U_m - U_L), \quad  D_1^+ = 0$.
        \item If $\displaystyle\sigma_- > 0$: 
        $\displaystyle D_1^- = 0, \quad D_1^+ = \sigma_- (U_m - U_L)$.
    \end{itemize}

    \item \textbf{1-Rarefaction ($\theta_L \ge \theta_m$)} bounded by speeds $\mu_-(U_L)$ and $\mu_-(U_m)$:
    \begin{itemize}
    \item Left-going ($\displaystyle  \mu_-(U_m) \le 0$): 
    \begin{eqnarray*} D_1^- &= & \mu_-(U_m) U_m - \mu_-(U_L) U_L - \int_{\mu_-(U_L)}^{\mu_-(U_m)} \widetilde{U}(\mu) \, d\mu, \\
        D_1^+ &= & 0 .
        \end{eqnarray*}
        \item Left-Right-going ($\displaystyle \mu_-(U_L) < 0 < \mu_-(U_m) $): 
        \begin{eqnarray*}
        D_1^- &=& -\mu_-(U_L) U_L - \int_{\mu_-(U_L)}^{0} \widetilde{U}(\mu) \, d\mu, \\
        D_1^+ &=& \mu_-(U_m) U_m - \int_{0}^{\mu_-(U_m)} \widetilde{U}(s) \, ds.          
        \end{eqnarray*}
        \item Right-going ($\displaystyle \mu_-(U_L) \ge 0$): 
        \begin{eqnarray*}
          D_1^- &=& 0, \\
        D_1^+ &=& \mu_-(U_m) U_m - \mu_-(U_L) U_L - \int_{\mu_-(U_L)}^{\mu_-(U_m)} \widetilde{U}(\mu) \, d\mu.  
        \end{eqnarray*}
    \end{itemize}
\end{itemize}

\bigskip

\noindent\textbf{2-Wave Contributions ($D_2^-, D_2^+$)}:
\begin{itemize}
    \item \textbf{2-Shock} ($\theta_R > \theta_m$) with speed $\displaystyle\sigma_+ = \frac{\rho_R \cos\theta_R - \rho_m \cos\theta_m}{\rho_R - \rho_m}$:
    \begin{itemize}
        \item If $\displaystyle \sigma_+ \le 0$: $\displaystyle D_2^- = \sigma_+ (U_R - U_m), \quad D_2^+ = 0$.
        \item If  $\displaystyle \sigma_+ > 0$: $\displaystyle D_2^- = 0, \quad D_2^+ = \sigma_+ (U_R - U_m)$.
    \end{itemize}

    \item \textbf{2-Rarefaction ($\theta_R \le \theta_m$)} bounded by speeds $\mu_+(U_m)$ and $\mu_+(U_R)$:
    \begin{itemize}
        \item Left-going ($\displaystyle \mu_+(U_R) \le 0$): 
        \begin{eqnarray*}
        D_2^- &=& \mu_+(U_R) U_R - \mu_+(U_m) U_m - \int_{\mu_+(U_m)}^{\mu_+(U_R)} \widetilde{U}(\mu) \, d\mu, \\
        D_2^+ &=& 0. 
        \end{eqnarray*}
        \item Left-Right-going  ($\displaystyle \mu_+(U_m) < 0 < \mu_+(U_R)$): 
        \begin{eqnarray*}
         D_2^- &=& -\mu_+(U_m) U_m - \int_{\mu_+(U_m)}^{0} \widetilde{U}(\mu) \, d\mu, \\
        D_2^+ &=& \mu_+(U_R) U_R - \int_{0}^{\mu_+(U_R)} \widetilde{U}(\mu) \, d\mu .   
        \end{eqnarray*}
        \item Right-going  ($\displaystyle\mu_+(U_m) \ge 0)$: 
        \begin{eqnarray*}
          D_2^- &=& 0, \\
        D_2^+ &=& \mu_+(U_R) U_R - \mu_+(U_m) U_m - \int_{\mu_+(U_m)}^{\mu_+(U_R)} \widetilde{U}(\mu) \, d\mu. 
        \end{eqnarray*}
    \end{itemize}
\end{itemize}

\begin{remark}
The integrals appearing at the contributions of 1-rarefaction waves of the form
$$
  \int_{\mu_-(U_L)}^{\mu_-(U_R)} \widetilde{U}(\mu) \, d\mu,$$
  where $U_L$ and $U_R$ belong to the same integral curve of the 1-characteristic field and $\mu_-(U_L) < \mu_-(U_R)$,
  can be written as follows
 \begin{equation}\label{eq:int_rar_1}
    \int_{\mu_-(U_L)}^{\mu_-(U_L)} \widetilde{U}(
    \mu) \, d\mu = \int_{\mu_-(U_L)}^{\mu_-(U_R)} \left[\begin{array}{c} \rho_L\exp(F_-(\chi_+(\mu))-F_-(\theta_L)) \\ \chi_+(\mu) \end{array} \right] \, d\mu,
\end{equation} 
where $F_-$ is given in Appendix \ref{appendix_2} and $\chi_+$ in \eqref{invav}.
Similarly, if $U_L$ and $U_R$ belong to the same integral curve of the 2-characteristic field and $\mu_+(U_L) < \mu_+(U_R)$ then
\begin{equation}\label{eq:int_rar_2}
    \int_{\mu_+(U_L)}^{\mu_+(U_R)} \widetilde{U}(\mu) \, d\mu = \int_{\mu_+(U_L)}^{\mu_+(U_R)} \left[\begin{array}{c} \rho_R\exp(F_+(\chi_-(\mu))-F_+(\theta_R)) \\ \chi_-(\mu) \end{array} \right] \, d\mu,
\end{equation}
where $F_+$  and $\chi_-$ are also given in Appendix \ref{appendix_2} and \eqref{invav}, respectively. Since an expression of these integrals in terms of elementary functions is not available, a quadrature formula is used here to approximate them.
\end{remark}

\section{High-order in-cell discontinuous reconstruction path-conservative methods for conservative systems in primitive variables}\label{sec:in_cell_method}

In order to numerically recover the right shock wave solutions of the SOH model (\ref{pde}), we need to enrich the Godunov method described in Section \ref{sec:Godunov}. In this section, we then define a high-order in time and space reconstruction operator which captures correctly shocks without numerical viscosity near them. For this to be achieved, the numerical scheme is given in the following semi-discrete form, similar to the path-conservative framework \cite{pares2006}:
%Similar to the path-conservative framework \cite{pares2006}, we write our numerical scheme in semi-discrete form:
\begin{equation}\label{eq:hopc}
    {U}_{i}'(t) = -\frac{1}{\Delta x}\left(\mathcal{D}_{i+\frac{1}{2}}^{-}(t) + \mathcal{D}_{i-\frac{1}{2}}^{+}(t) + \dashint_{x_{i-\frac{1}{2}}}^{x_{i+\frac{1}{2}}} \mathcal{A}(P^n_{i}(x,t))\frac{\partial}{\partial x} P^n_{i}(x,t)dx\right),\quad t\in (t_n, t_{n+1}],
\end{equation}
where:
\begin{itemize}
    % \item $\Delta x$ and $\Delta t$ are the space and the time steps respectively. They are assumed to be constant in the description of the numerical method for simplicity.
    % \item $I_{i} = \left[x_{i-\frac{1}{2}}, x_{i-\frac{1}{2}}\right]$ are the computational cells with lenght $\Delta x$.
    % \item $t_n = n\Delta t$, $n=0,1,2,...$.
	\item $U_{i}(t) \approx \displaystyle \frac{1}{\Delta x}\int_{x_{i+\frac{1}{2}}}^{x_{i-\frac{1}{2}}} U(x,t) \,dx$ is the approximation to the cell average of the solution at the $i$-th cell in time $t$;
    \item $P^n_{i}(x,t)$ is a space-time reconstruction operator, which may be discontinuous, to be defined in Subsection \ref{sec:reconstruction_operator} so that numerical viscosity is removed close to shocks in order to capture them correctly;
	\item $\mathcal{D}^\pm_{i+1/2}(t) = {D}^\pm(P^n_{i}(x_{i+1/2}, t),P^n_{i+1}(x_{i+1/2}, t))$, where ${D}^\pm$  are the fluctuations of the Godunov method defined in Section \ref{sec:Godunov}.

\end{itemize}

Once the numerical approximations $U_i^n$ at time $t^n = n \Delta t$ have been computed, the ODE system \eqref{eq:hopc} has to be solved with an adequate time stepping to update the numerical solution. 

This section is decomposed as follows: the reconstruction operator $P_i^n$ is defined in Subsection \ref{sec:reconstruction_operator}. Subsection \ref{sec:markingcells} presents the marking of cells. Then the fully-discrete form of (\ref{eq:hopc}) is given in Subsection \ref{sec:fullydiscrete}. Afterwards, the time step is described in Subsection \ref{sub:time_step}. The numerical schemes are then given at order 1 in Subsection \ref{sec:1storder}, and order 2 in Subsection \ref{sec:2ndorder}. Finally, the theorem that states the shock-capturing property of the numerical method is recalled in Subsection \ref{sec:conv}.

\subsection{Definition of reconstruction operators $P^n_{i}(x,t)$}\label{sec:reconstruction_operator}
The first step is to identify the cells $I_i$ where a shock wave is present or may be generated. More precisely, we mark those cells for which the solution of the
Riemann problem associated with \eqref{pde}, with initial data
\begin{equation}\label{eq:Riemann}
U(x,0) =
\begin{cases}
U_{i-i}^n, & x < 0, \\
U_{i+1}^n, & x > 0,
\end{cases}
\end{equation}
involves only one shock wave. Since system \eqref{pde} have only two characteristics fields, we have only four possibilities: 
\begin{itemize}
    \item Isolated 1-shock.
    \item Isolated 2-shock.
    \item 1-shock + 2-rarefaction.
    \item 1-rarefaction + 2-shock.
\end{itemize}

In the rest of the cases we have isolated rarefactions, two shocks or two rarefactions. In these cases we do not mark the cell.

We denote by $\mathcal{M}_n$ the set of indices corresponding to the marked cells, namely
\begin{equation}
\mathcal{M}_n = \left\{\, i \ \  \text{s.t. the solution of the Riemann problem \eqref{pde},
\eqref{eq:Riemann} contains only one shock wave} \,\right\}.
\end{equation}

The second step is to select a space-time reconstruction operator that, given the cell averages $U^n_i$ at the cells of a smooth function $U(x,t)$ at time $t_n$, provides smooth approximations of $U$ at the cells 
$$
U(x,t) \approx Q_i(x,t) = Q_i(x,t; \{ U^n_j \}_{j \in \mathcal{S}_{i} }), \quad x \in [x_{i-1/2}, x_{i+1/2}], \ t \in [t_n, t_{n+1}],
$$
where $\mathcal{S}_{i}$ is the stencil, i.e. the set of indices of the cells whose averages are used to compute $Q_i$. In the case of first order methods, the piecewise constant reconstruction operator will be considered. For higher order methods,  Taylor expansions and the Cauchy-Kowaleskaya procedure will be followed to compute these reconstruction operators as in the ADER approach \cite{titarev2005ader, dumbser2009ader}.

On the basis of the marked cells and the standard reconstruction operator, $P_{i}^{n}(x,t)$ is then defined as follows:
\begin{itemize}

\item If $i \in \mathcal{M}_n$ then the reconstruction operator is defined by
        \begin{equation}\label{eq:in_cell_reconstruction}
        P_{i}^{n}(x,t) = 
\begin{cases}
U^n_{i,l} & \text{ if $x \leq x_{i - 1/2} + d_{i}^{n} \Delta x + \sigma_i^n (t -t^n) $,}\\
U^n_{i,r} & \text{otherwise},
\end{cases}
\end{equation}
where 
\begin{itemize}

\item  $U^n_{i,l}$, $U^n_{i,r}$ are two states and $\sigma_i^n$ a scalar to be chosen so that if $U^n_{i-1}$ and
$U^n_{i+1}$ can be  linked by an admissible discontinuity with speed $\sigma$, then
\begin{equation}\label{speedstates}
U^n_{i,l} = U^n_{i-1}, \quad U^n_{i,r} = U^n_{i+1}, \quad \sigma_i^n = \sigma.
\end{equation}
See Figure \ref{fig:In_cell_discontinuous_reconstruction_operator} where the case of a 2-shock reconstruction is illustrated.

\item $d_i^n$ is chosen so that
\begin{equation}\label{conservation}
d_i^n u_{i,l, k}^n + (1 -  d_i^n) u_{i,r, k}^n = u_{i, k}^n, 
\end{equation}
where $u_{i,k}^n$ represents the $k$th component of $U_i^n$, for some index $k \in \{ 1, 2 \}$. In practice we use the $\rho$ component.

\end{itemize}
Observe that this in-cell discontinuous reconstruction can only be done
if $0 \leq d_i^n \leq 1$, i.e. if
\begin{equation}\label{eq:d_property}
    0 \leq \frac{u_{i,r,k}^n - u_{i,k}^n}{u_{i,r,k}^n - u_{i,l,k}^n} \leq 1,
\end{equation}
otherwise the index $i$ is removed from the set $\mathcal{M}_n$.
\item If $i \notin \mathcal{M}_n$ and $\mathcal{S}_{i}\cap  \mathcal{M}_n \neq \emptyset$, then 
    \begin{equation}\label{eq:Constant_reconstruction_operator}
        P_{i}^{n}(x,t)= U_{i}^{n}.
    \end{equation}
See Figure \ref{fig:Constant_reconstruction_operator}.
\item Otherwise (i.e.,  if $\mathcal{S}_{i}\cap  \mathcal{M}_{n} = \emptyset$)
$$
P_i^n(x,t) = Q_i^n(x,t; \{ U_j^n \}_{j \in \mathcal{S}_i} ).
$$

\begin{figure}[htpb]
		\begin{subfigure}{0.5\textwidth}
			\includegraphics[width=0.9\linewidth]{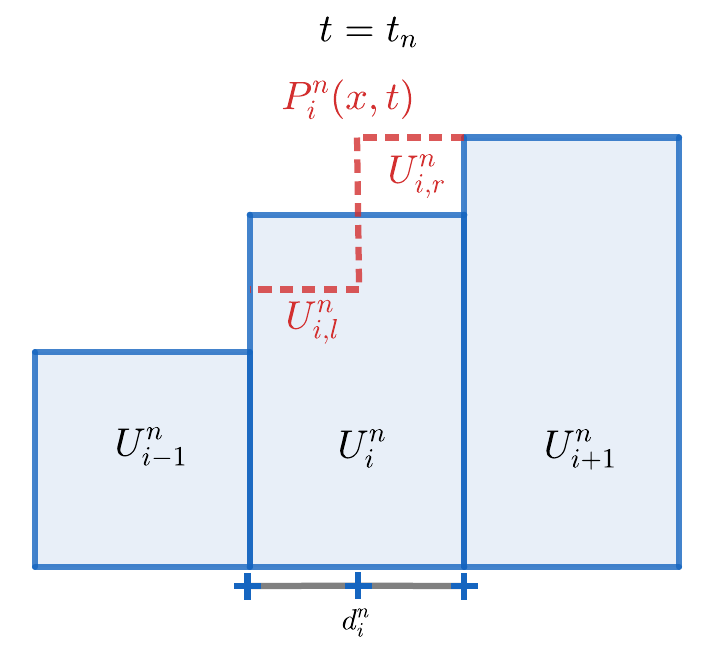}
			\caption*{In-cell discontinuous reconstruction operator at time $t = t^{n}$.}
		\end{subfigure}
		\begin{subfigure}{0.5\textwidth}
				\includegraphics[width=0.9\linewidth]{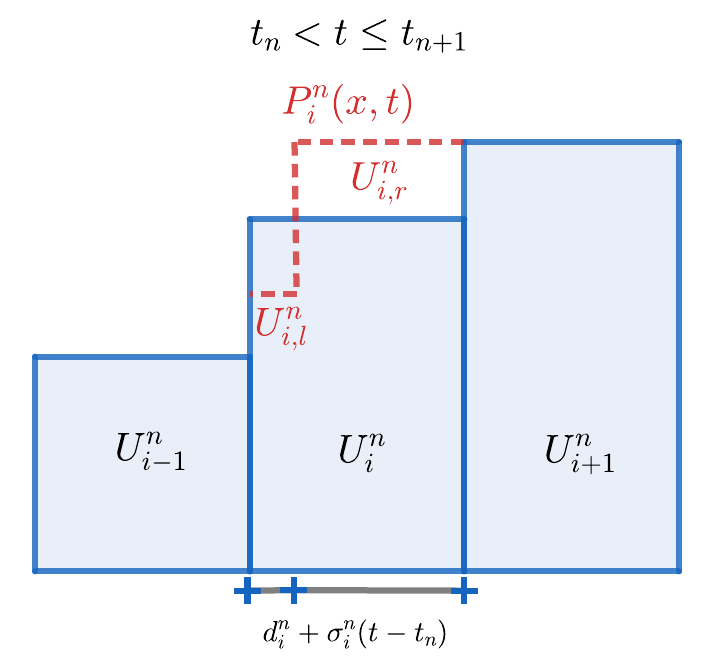}
				\caption*{In-cell discontinuous reconstruction operator at intermediate time $t^{n} < t < t^{n+1}$}
		\end{subfigure}
		\caption{In-cell discontinuous reconstruction operator.}
		\label{fig:In_cell_discontinuous_reconstruction_operator}
	\end{figure}
        \begin{figure}[htpb]
        \centering
			\includegraphics[width=0.5\linewidth]{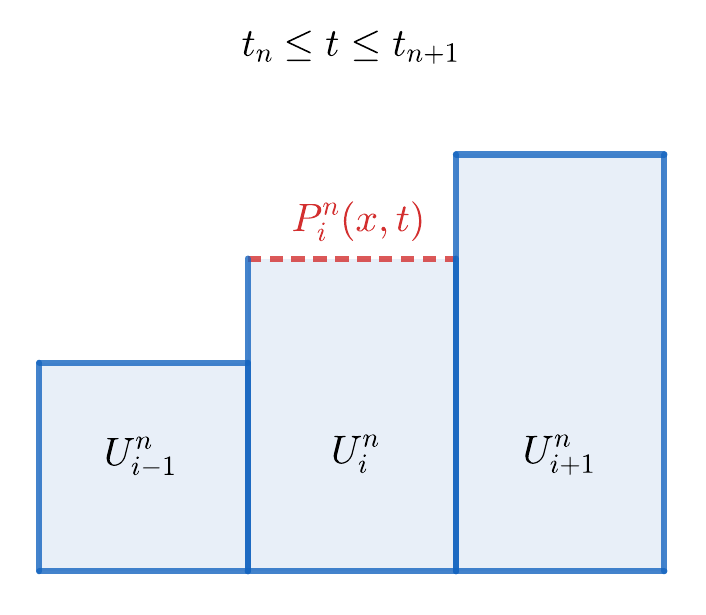}
			\caption{Constant reconstruction operator.}
		\label{fig:Constant_reconstruction_operator}
	\end{figure}
\end{itemize}

\begin{remark} In practice, in the case $i \in \mathcal{M}_n$, the index $k$th related to the mass-conservation equation is chosen.

\end{remark}

\subsection{How to mark a cell? Choices of $\sigma_{i}^{n}, U_{i,l}^{n}$ and $U_{i,r}^{n}$} \label{sec:markingcells}

In \cite{Vicsek_Riemann}, the solutions of the Riemann problem associated with \eqref{pde} were obtained. This information is used to choose the states and the speed of the discontinuity in \eqref{eq:in_cell_reconstruction}. Let us consider again the Riemann problem \eqref{eq:Riemann} whose intermediate state is given by $W_{i,m}^n = [\rho_{i,m}^{n}, \theta_{i,m}^n]^T$. Then for each cell $I_i$ of the mesh
\begin{enumerate}

\item If $\theta_m > \theta_L$ and $ \theta_m \geq \theta_R$ (1-shock + 2-rarefaction or isolated 1-shock wave), then the reconstruction operator is given by \eqref{eq:in_cell_reconstruction} with
\begin{equation}\label{eq:reconstructed_states_and_speed_1}
   U^n_{i,l} = U^n_{i-1}, \quad U^n_{i,r} = U_{i,m}^n = [\rho_{i,m}^n,\theta_{i,m}^n]^T , \quad   \sigma_i^n = \sigma_{i,m}^{n} = \frac{\rho_{i,m}^n \cos\theta_{i,m}^n - \rho_{i-1}^{n} \cos\theta_{i-1}^{n}}{\rho_{i,m}^{n} - \rho_{i-1}^{n}}.
\end{equation} 
\item If $\theta_m \leq \theta_L$ and $ \theta_m < \theta_R$ (1-rarefaction and a 2-shock wave or isolated 2-shock wave), then the reconstruction operator is given by \eqref{eq:in_cell_reconstruction} with
\begin{equation}\label{eq:reconstructed_states_and_speed_2}
   U^n_{i,l} = U_{i,m}^n = [\rho_{i,m}^n,\theta_{i,m}^n]^T, \quad U^n_{i,r} = U_{i+1}^n, \quad   \sigma_i^n = \sigma_{i,m}^{n} = \frac{\rho_{i+1}^n \cos\theta_{i+1}^n - \rho_{i,m}^{n} \cos\theta_{i,m}^{n}}{\rho_{i+1}^{n} - \rho_{i,m}^{n}}.
\end{equation} 
\item Otherwise, no shock wave or two shock waves are located in the cell $I_i$, so the constant reconstruction operator \eqref{eq:Constant_reconstruction_operator} is considered.
\end{enumerate}
Summing up, the set $\mathcal{M}_n$ is given by
$$
\mathcal{M}_{n} = \{ i \ \text{ s.t. the Riemann problem \eqref{eq:Riemann} satisfies case 1 or case 2}\}.
$$

\begin{remark}
    In practice, for cases 1 \eqref{eq:reconstructed_states_and_speed_1} and 2 \eqref{eq:reconstructed_states_and_speed_2} we use the following conditions, respectively:
    \begin{enumerate}
        \item $\theta_m > \theta_L +\epsilon$ and $ \theta_m \geq \theta_R - \epsilon$.
        \item $\theta_m \leq \theta_L - \epsilon$ and $ \theta_m +\epsilon < \theta_R$,
    \end{enumerate}
    where we use $\epsilon = 0.01$. This allows us to mark the cell near the formation of a shock.
\end{remark}

\subsection{Time step}\label{sub:time_step}
Following \cite{pimentelgarcia2022cell, pimentel2024drmood} the time step $\Delta t_n$ is chosen as follows:
\begin{equation}\label{restricted_time_step}
    \Delta t_n = \min(\Delta t^c_n, \Delta t^r_n).
\end{equation}
Here
\begin{equation}
    \Delta t_n^c = CFL\cdot \left( \frac{ \Delta x}{\max_i\left\{|\mu_{-}(\theta_i^n)|, |\mu_{+}(\theta_i^n)|\right\}} \right),
\end{equation}
where $CFL \in (0,1)$ is the stability parameter, $\mu_{\pm}(\theta_i^n)$ represent the eigenvalues of
$\mathcal{A}(U_i^n)$ given by (\ref{av1}), and
\begin{equation}
    \Delta t_n^r = \min_{i \in \mathcal{M}_n} \begin{cases}
    \displaystyle \frac{1 - d_i^n}{|\sigma_i^n|}\Delta x, & \quad if \quad \sigma_i^n>0, \\ \\
    \displaystyle \frac{d_i^n}{|\sigma_i^n|}\Delta x, & \quad if \quad \sigma_i^n<0,
    \end{cases} 
\end{equation}
where $d_{i}^{n}$ and $\sigma_{i}^{n}$
are given in \eqref{conservation} and \eqref{eq:reconstructed_states_and_speed_1}-\eqref{eq:reconstructed_states_and_speed_2}, respectively, only for marked cells. Observe that, besides the stability requirement, this choice of time step ensures that no discontinuous reconstruction leaves a marked cell. Note that this time restriction could be avoided as in \cite{chalons2020path} but, for simplicity, we keep it.

\subsection{Fully discrete form} \label{sec:fullydiscrete}

$U_i^{n+1}$ is computed by integrating \eqref{eq:hopc} in the interval $t \in [t_{n}, t_{n+1}]$: 
\begin{equation}
    {U}_{i}^{n+1} = U_{i}^{n}-\int_{t_{n}}^{t_{n+1}}\frac{1}{\Delta x}\left(\mathcal{D}_{i+\frac{1}{2}}^{-}(t) + \mathcal{D}_{i-\frac{1}{2}}^{+}(t) + \dashint_{x_{i-\frac{1}{2}}}^{x_{i+\frac{1}{2}}} \mathcal{A}(P_{i}^{n}(x,t))\frac{\partial}{\partial x} P_{i}^{n}(x,t)dx\right)dt,
\end{equation}
where $t_{n+1}=t_n + \Delta t_n$ was defined in Subsection \ref{sub:time_step}.  Two quadrature formulas
$$
\int_{x_{i-1/2}}^{x_{i+1/2}} f(x) \, dx \approx \Delta x \sum_{m = 1}^{M_{x}} \beta_{m}^{x} f(x_i^{m}), \quad
\int_{t^{n}}^{t^{n+1}} g(t) \, dt \approx \Delta t_n \sum_{l = 1}^{M_{t}} \beta_{l}^{t} g(t^{n,l}),
$$
are used to approximate the integrals in space and time, respectively. Taking into account the expression of the reconstruction operators we obtain:

\begin{equation}\label{eq:discrete1}
    U_{i}^{n+1} = \begin{cases}
        \displaystyle U_{i}^{n} - \frac{\Delta t_n}{\Delta x}\left(\mathcal{D}_{i+\frac{1}{2}}^{-} + \mathcal{D}_{i-\frac{1}{2}}^{+} + \sigma_{i}^{n}(U_{i,r}^{n}- U_{i,r}^{n}) \right) & \text{if $i \in \mathcal{M}_n$,}\\
        \displaystyle U_{i}^{n} - \frac{\Delta t_n}{\Delta x}\left(\mathcal{D}_{i+\frac{1}{2}}^{-} + \mathcal{D}_{i-\frac{1}{2}}^{+}\right) & \text{if $i \notin \mathcal{M}_n$ and  $\mathcal{S}_{i}\cap \mathcal{M}_n \neq \emptyset$,}\\
        \displaystyle U_{i}^{n} - \frac{\Delta t_n}{\Delta x}\sum_{l=1}^{M_{t}} \beta_{l}^{t}\left(\mathcal{D}_{i+\frac{1}{2}}^{-}(t^{n,l}) + \mathcal{D}_{i-\frac{1}{2}}^{+}(t^{n,l})\right) & \\
        \qquad \displaystyle- \Delta t_n \sum_{l=1}^{M_{t}} \sum_{m=1}^{M_{x}} \beta_{l}^{t} \beta_{m}^{x}\mathcal{A}(P_{i}^{n}(x_{i}^{m},t^{n,l}))
    \partial_x P_{j}^{n}(x_{i}^{m},t^{n,l}) & \text{if $\mathcal{S}_{i} \cap \mathcal{M}_{n} = \emptyset$.} 
    \end{cases}
\end{equation}

\begin{remark}
    Observe that this method is basically a standard high-order Godunov method in the smooth region and the in-cell discontinuous reconstruction developed in \cite{chalons2020path, pimentelgarcia2022cell, pimentel2024drmood} in the presence of shocks for pure nonconservative systems.
\end{remark}

\subsection{First-order method} \label{sec:1storder}

In this case of first-order methods, the space-time reconstruction operator is the piecewise constant one, i.e.
$$Q_{i}^{n}(x,t) = U_{i}^{n},$$
so that the first-order in-cell discontinuous reconstruction path-conservative method \eqref{eq:discrete1} method is reduced to:
\begin{equation}\label{eq:discretefo}
    U_{i}^{n+1} = \begin{cases}
        \displaystyle U_{i}^{n} - \frac{\Delta t_n}{\Delta x}\left(\mathcal{D}_{i+\frac{1}{2}}^{-} + \mathcal{D}_{i-\frac{1}{2}}^{+}+ \sigma_{i}^{n}(U_{i,r}^{n}- U_{i,r}^{n}) \right) & \text{if $i \in \mathcal{M}_n$,}\\
        \displaystyle U_{i}^{n} - \frac{\Delta t_n}{\Delta x}\left(\mathcal{D}_{i+\frac{1}{2}}^{-} + \mathcal{D}_{i-\frac{1}{2}}^{+}\right) & \text{if $i \notin \mathcal{M}_n$.}\\
    \end{cases}
\end{equation}

\subsection{Second-order method} \label{sec:2ndorder}

The MUSCL-Hancock space-time reconstruction operator (see\cite{van1974towards, van1984relation}) is considered in this case:
\begin{equation}
Q_{i}^{n}(x,t)= U_{i}^{n} + \widetilde{\partial_xU}_{i}^{n}(x-x_{i})- \mathcal{A}(U_{i}^{n})\widetilde{\partial_xU}_{i}^{n}(t-t^{n}),
    \end{equation}
and the stencil is $\mathcal{S}_i = \{i-1,i,i+1\}$.
Here, $\widetilde{\partial_xU}_{j}^{n}$ is the \textit{minmod} approximation of the first order spacial derivative of 
$U$ at $x_i$ at time $t_n$, whose $k$th component is given by
$$
\left(\widetilde{\partial_xU}_{i}^{n}\right)_k = {\rm{minmod}}\left( \frac{u_{i+1,k}^{n}-u_{i,k}^{n}}{\Delta x}, \frac{u_{i,k}^{n}-u_{i-1,k}^{n}}{\Delta x}\right),$$
where $u_{i,k}^n$ represents again the $k$th component of $U_i^n$ and 
$$\rm{minmod}(a,b) = \begin{cases}
     \min\{a,b\} & \text{if $a,b>0$,}  \\
     \max\{a,b\} & \text{if $a,b<0$,}  \\
     0 & \text{otherwise.}
     \end{cases}$$
The mid-point rule is chosen as quadrature formula both in time and space, so that the second-order in-cell discontinuous reconstruction path-conservative method \eqref{eq:discrete1} is given by:
\begin{equation}\label{eq:discrete2}
    U_{i}^{n+1} = \begin{cases}
        \displaystyle U_{i}^{n} - \frac{\Delta t_n}{\Delta x}\left(\mathcal{D}_{i+\frac{1}{2}}^{-} + \mathcal{D}_{i-\frac{1}{2}}^{+}+ \sigma_{i}^{n}(U_{i,r}^{n}- U_{i,r}^{n}) \right) & \text{if $i \in \mathcal{M}_n$,}\\
        \displaystyle U_{i}^{n} - \frac{\Delta t_n}{\Delta x}\left(\mathcal{D}_{i+\frac{1}{2}}^{-} + \mathcal{D}_{i-\frac{1}{2}}^{+}\right) & \text{if $i \notin \mathcal{M}_n$ and $\mathcal{S}_{i} \cap \mathcal{M}_{n} \neq \emptyset$,}\\
        \displaystyle U_{i}^{n} - \frac{\Delta t_n}{\Delta x}\left(\mathcal{D}_{i+\frac{1}{2}}^{-}(t_{n+\frac{1}{2}}) + \mathcal{D}_{i-\frac{1}{2}}^{+}(t_{n+\frac{1}{2}}) + \Delta x \mathcal{A}(U_{i}^{n+\frac{1}{2}})\widetilde{\partial_xU}_{i}^{n}\right) & \text{if $\mathcal{S}_{i} \cap \mathcal{M}_{n} = \emptyset$,} 
    \end{cases}
\end{equation}
where $\displaystyle t_{n+\frac{1}{2}} = t_n + \frac{\Delta t_n}{2}$ so
    $$\mathcal{D}_{i+\frac{1}{2}}^{\pm}(t_{n+\frac{1}{2}}) = \mathcal{D}^{\pm}\left(P_{i}^{n}\left(x_{i+\frac{1}{2}},t_n + \frac{\Delta t_n}{2}\right),P_{i+1}^{n}\left(x_{i+\frac{1}{2}},t_n + \frac{\Delta t_n}{2}\right)\right) $$
    and
    $$U_{i}^{n+\frac{1}{2}} =  P_{i}^{n}\left(x_{i},t_n + \frac{\Delta t_n}{2}\right).$$

\subsection{Shock-capturing and Rankine-Hugoniot compatibility} \label{sec:conv}

We recall the theorem in \cite{chalons2020path, pimentelgarcia2022cell, pimentel2024drmood} that states that the scheme \eqref{eq:discrete1} captures exactly isolated shocks in the following way:

\begin{theorem}\label{th}
Assume that $U_l$ and $U_r$ can be joined by an entropy shock of speed $\sigma$ satisfying \eqref{RH}. Then, the numerical method  \eqref{eq:discrete1} provides an exact numerical solution
of the Riemann problem with initial conditions
$$
U(x,0) = \begin{cases}
U_l & \text{if $x < 0$,} \\
U_r & \text{otherwise,}
\end{cases} 
$$
in the sense that
\begin{equation}\label{presave}
U_i^n = \frac{1}{\Delta x}\int_{x_{i-1/2}}^{x_{i+1/2}} U(x,t_n) \, dx, \quad \forall i, n,
\end{equation}
where $U(x,t)$ is the exact solution.

\end{theorem}

\section{Numerical tests}\label{sec:numerical_tests}

This section is devoted to show the convergence to the exact solution when using the Godunov in-cell discontinuous reconstruction method and to compare it with the Godunov method without the in-cell discontinuous reconstruction. The notation is going to be:
\begin{itemize}
    \item O1\_God\_DR: First-order Godunov scheme with in-cell discontinuous reconstruction.
    \item O1\_God\_noDR: First-order Godunov scheme without in-cell discontinuous reconstruction.
    \item O2\_God\_DR: Second-order Godunov scheme with in-cell discontinuous reconstruction.
    \item O2\_God\_noDR: Second-order Godunov scheme without in-cell discontinuous reconstruction.
    \item O1\_HLL\_Viscous: First-order HLL-type viscous shock-capturing scheme developed in \cite{compain:hal-05490904}. This method consists on a Godunov-type discretization for the SOH model (\ref{pde}) with a numerical viscous correction. The numerical viscous parameters used for the following numerical simulations are the ones prescribed in \cite{compain:hal-05490904}.
\end{itemize}
In all tests a $CFL=0.5$ is used together with free boundary conditions and the spatial domain $I=[0,10]$. We are going to consider Riemann problems, so the initial conditions will be of the form:
\begin{equation}
    U_{0}(x) = \begin{cases}
        U_{L} = (\rho_L, \theta_L)^{T} & \text{if } x\leq 5, \\
        U_{R} = (\rho_R, \theta_R)^{T} & \text{if } x> 5,
    \end{cases}
\end{equation}
where the values $\rho_L, \theta_L, \rho_R, \theta_R$ are given in Table \ref{tab:initial}.

\begin{table}[!h]
\centering

\begin{tabular}{|c|c|c|c|c|c|c|c|}
\hline
\multicolumn{8}{|c|}{Tests taken from \cite{compain:hal-05490904}} \\
\hline
Tests &
$\rho_L$ &
$\theta_L$ &
$\rho_R$ &
$\theta_R$ & N & $c$ & $\lambda$\\
\hline
Test 1 & 1 & 0.5 & 3.078009 & 1.281671 & 1000 & 0.7779710678272852 & 0.2238681162728627\\
Test 2 & 1 & 0.5 & 1.867375 & 0.915616 & 200 & 0.7779710678272852 & 0.2238681162728627\\
Test 3 & 1 & 0.5 & 4.445344 & 1.371887 & 1000 & 0.7779710678272852 & 0.2238681162728627\\
\hline
\multicolumn{8}{|c|}{Tests taken from \cite{Vicsek_Riemann}} \\
\hline
Test 4 & 1 & 1 & 1 & 2 & 1000 & 0.8 & 0.2\\
Test 5 & 1 & 2 & 0.5 & 1 & 200 & 0.8 & 0.2\\
Test 6 & 0.5 & 1.5 & 5 & 2 & 1000 & 0.8 & 0.2\\
\hline
\end{tabular}
\caption{Initial conditions for the test cases }
\label{tab:initial}
\end{table}

\subsection{Tests taken from \cite{compain:hal-05490904}}

The justification of the election of the parameters $c$ and $\lambda$ can be found in \cite{Degond2012, Degond2008, motsch2011numerical}.

\paragraph{Test 1: Stationary isolated 1-shock}

This test consists on an isolated 1-shock with speed $\sigma_1 = 0$. Figure \ref{fig:Test1} compares the Godunov methods with and without the in-cell discontinuous reconstruction (DR), together with the HLL viscous scheme. The schemes without DR are not able to preserve the stationary shock: the discontinuity moves away from its initial position, an incorrect intermediate state appears and travels to the right, and spurious oscillations are produced near the shock. In contrast, both schemes with DR preserve the stationary shock exactly. The HLL viscous scheme keeps the shock at the correct position, but smears it over a few cells and generates a small spurious oscillation, in a region where the exact solution is constant.

\begin{figure}[htpb]
	\begin{subfigure}[h]{0.5\textwidth}
		\centering
			\includegraphics[width=1.1\linewidth]{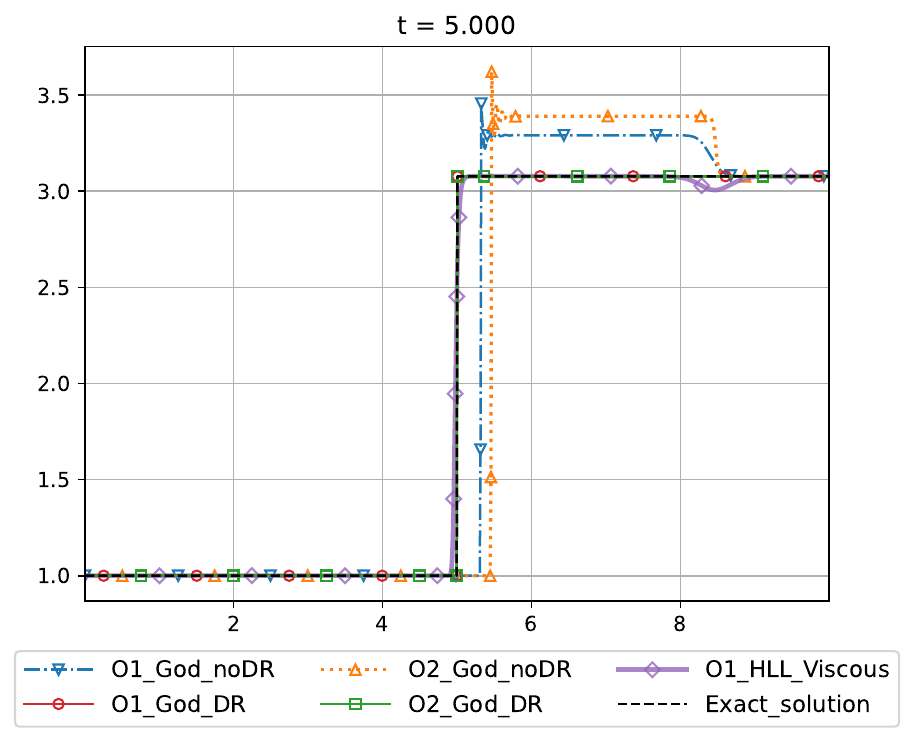}
            \caption{Variable $\rho$}
	\end{subfigure}
	\begin{subfigure}[h]{0.5\textwidth}
		\centering
			\includegraphics[width=1.1\linewidth]{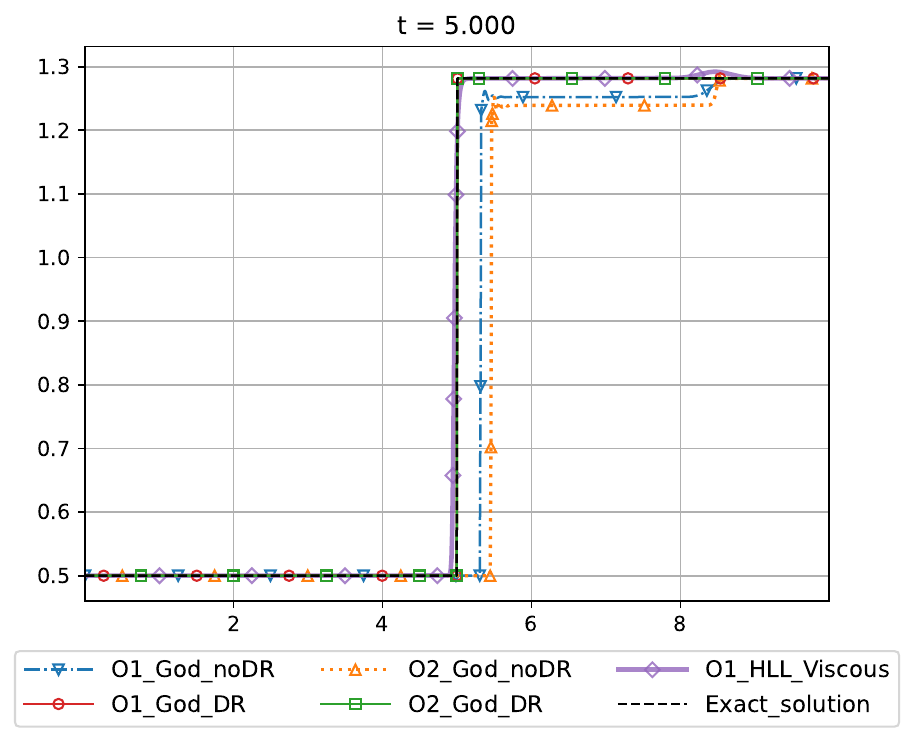}
            \caption{Variable $\theta$}
	\end{subfigure}
	\caption{Numerical results for test 1: stationary isolated 1-shock.}
	\label{fig:Test1}
\end{figure}

\paragraph{Test 2: isolated 1-shock with positive speed}

This test consists on isolated 1-shock with speed $\sigma_1 = 0.3$. In Figure \ref{fig:Test2}, the Godunov schemes without DR give a shock with the wrong speed and an incorrect post-shock state. With DR, both schemes capture the shock at the right position and without numerical diffusion as expected from Theorem \ref{th}. The HLL viscous scheme gives the correct shock location, but the discontinuity is strongly smeared again and the post-shock state is not accurately recovered.

\begin{figure}[htpb]
	\begin{subfigure}[h]{0.5\textwidth}
		\centering
			\includegraphics[width=1.1\linewidth]{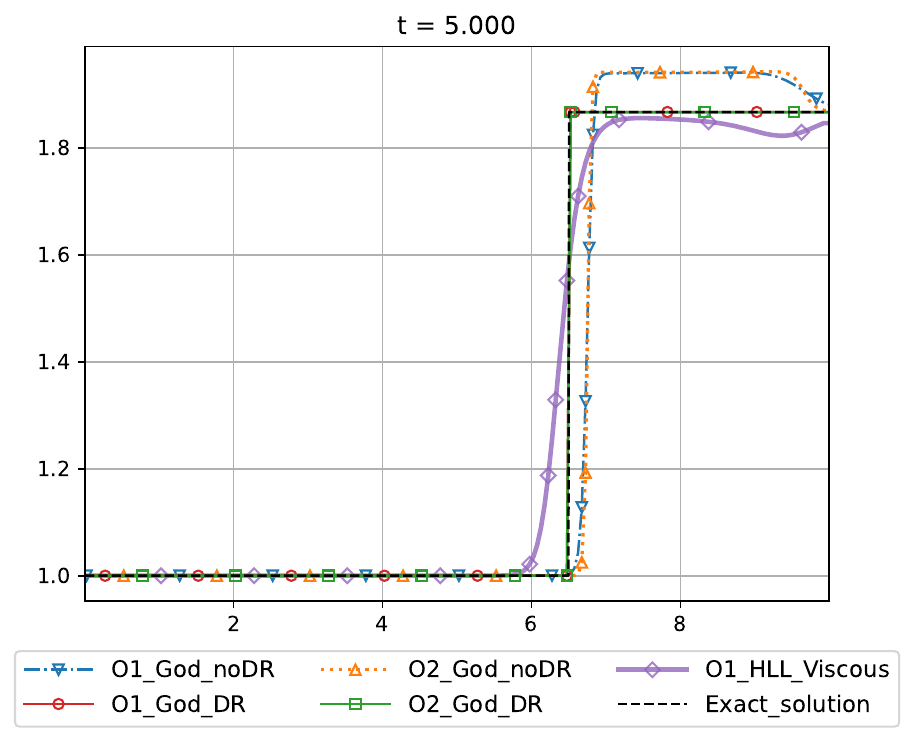}
            \caption{Variable $\rho$}
	\end{subfigure}
	\begin{subfigure}[h]{0.5\textwidth}
		\centering
			\includegraphics[width=1.1\linewidth]{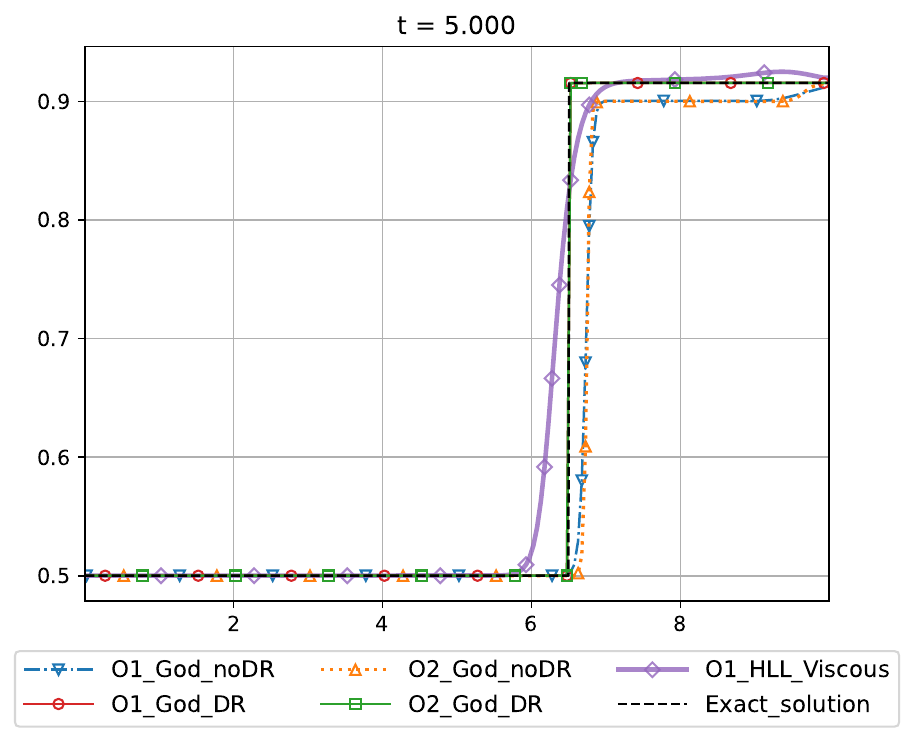}
            \caption{Variable $\theta$}
	\end{subfigure}
	\caption{Numerical results for test 2: isolated 1-shock with positive speed.}
	\label{fig:Test2}
\end{figure}

\paragraph{Test 3: 1-shock + 2-rarefaction}

This test consists on a 1-shock with speed $\sigma_1 = -0.8665350302618553$ followed by a 2-rarefaction linked through an intermediate state \\ $U_{m} = (\rho_m, \theta_m)^{T} = (2.64567223988045, 2.3314594288198824)^{T}$. Figure \ref{fig:Test3} shows that the Godunov schemes without DR completely miss the intermediate state and produce a single wrong discontinuity instead of the shock-rarefaction structure. With DR, both the shock and the intermediate state are correctly captured, and the second-order scheme resolves the rarefaction more sharply than the first-order one. The HLL viscous scheme captures the shock at the right position, although it is smeared, and it introduces a small perturbation in the intermediate state. Its rarefaction is about as diffusive as with the first-order DR scheme. In Table \ref{tab:Test3_convergence} and Figure \ref{fig:Test3_convergence} we observe the errors obtained with O1\_God\_DR and O2\_God\_DR for different meshes. Since the solution contains a shock, Table \ref{tab:Test3_convergence} and Figure \ref{fig:Test3_convergence} are not intended to measure the order of the schemes, but to show numerically that the methods with in-cell discontinuous reconstruction converge to the exact solution, even though no theoretical result guarantees it yet.

\begin{figure}[htpb]
	\begin{subfigure}[h]{0.5\textwidth}
		\centering
			\includegraphics[width=1.1\linewidth]{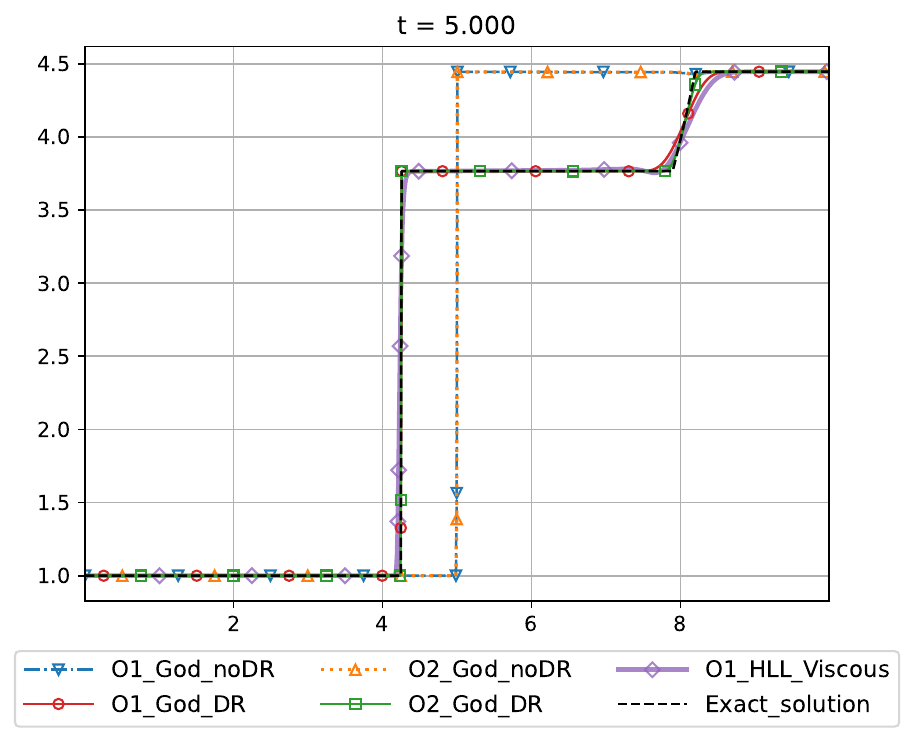}
            \caption{Variable $\rho$}
	\end{subfigure}
	\begin{subfigure}[h]{0.5\textwidth}
		\centering
			\includegraphics[width=1.1\linewidth]{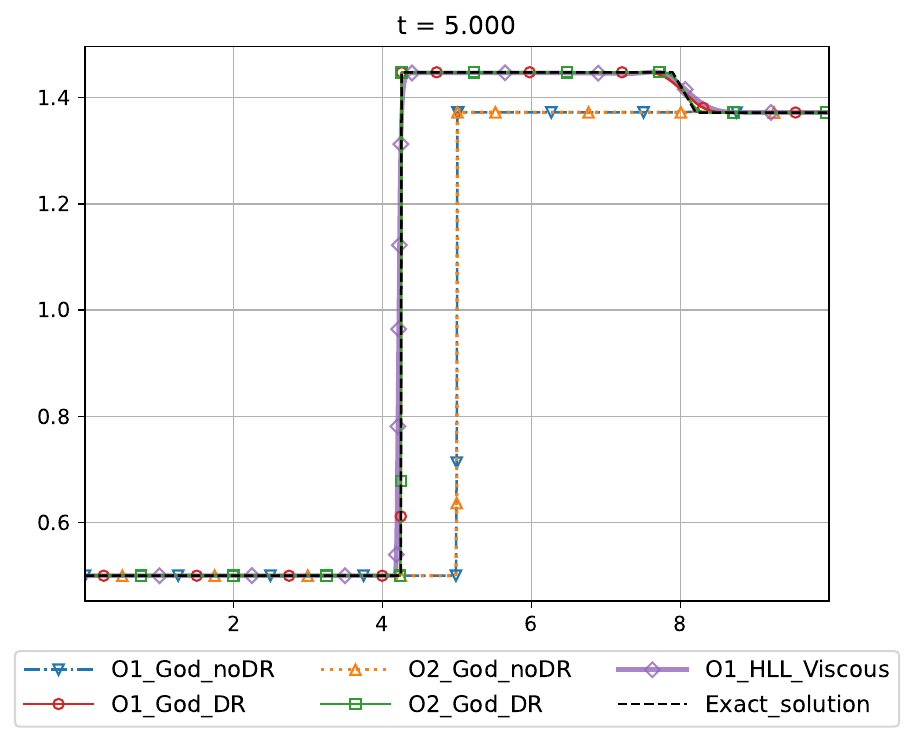}
            \caption{Variable $\theta$}
	\end{subfigure}
	\caption{Numerical results for test 3: 1-shock + 2-rarefaction.}
	\label{fig:Test3}
\end{figure}

\begin{table}[!h]
\centering

\begin{tabular}{|c|c|c|c|c|c|c|c|c|}
\hline
 & \multicolumn{4}{c|}{O1\_God\_DR} & \multicolumn{4}{c|}{O2\_God\_DR} \\
\hline
N &
$L^1$ error $\rho$ & Order &
$L^1$ error $\theta$ & Order &
$L^1$ error $\rho$ & Order &
$L^1$ error $\theta$ & Order \\
\hline
100  & $3.04\cdot 10^{-1}$ & --   & $5.28\cdot 10^{-2}$ & --   & $1.51\cdot 10^{-1}$ & --   & $2.77\cdot 10^{-2}$ & --   \\
200  & $1.97\cdot 10^{-1}$ & 0.62 & $3.18\cdot 10^{-2}$ & 0.73 & $8.54\cdot 10^{-2}$ & 0.82 & $1.51\cdot 10^{-2}$ & 0.87 \\
400  & $1.28\cdot 10^{-1}$ & 0.63 & $1.94\cdot 10^{-2}$ & 0.71 & $4.62\cdot 10^{-2}$ & 0.89 & $8.04\cdot 10^{-3}$ & 0.91 \\
800  & $8.29\cdot 10^{-2}$ & 0.62 & $1.20\cdot 10^{-2}$ & 0.69 & $2.42\cdot 10^{-2}$ & 0.93 & $4.18\cdot 10^{-3}$ & 0.94 \\
1600 & $5.40\cdot 10^{-2}$ & 0.62 & $7.60\cdot 10^{-3}$ & 0.66 & $1.25\cdot 10^{-2}$ & 0.95 & $2.17\cdot 10^{-3}$ & 0.95 \\
\hline
\end{tabular}
\caption{Test 3: $L^1$ errors and convergence orders at $t=5$ of the first- and second-order Godunov schemes with in-cell discontinuous reconstruction.}
\label{tab:Test3_convergence}
\end{table}

\begin{figure}[htpb]
		\centering
\includegraphics[width=1\linewidth]{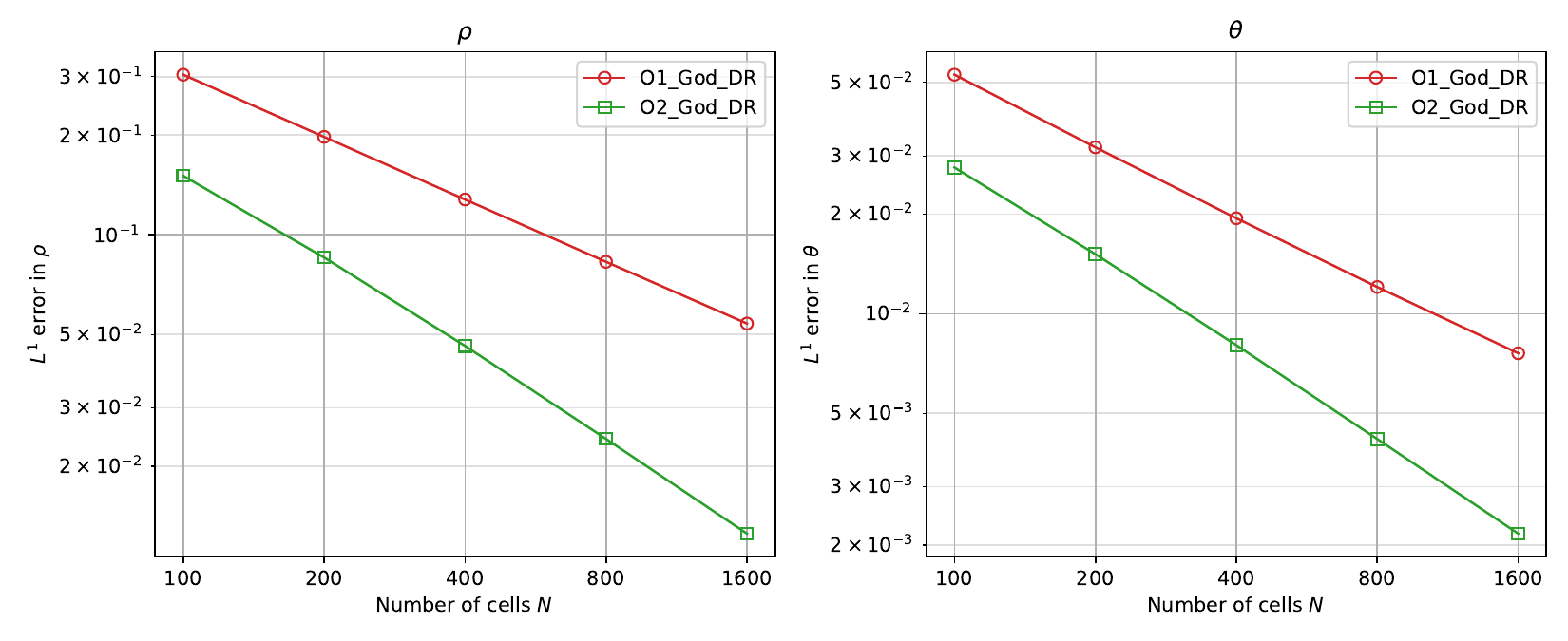}
	\caption{Test 3: $L^1$ errors vs number of cells.}
	\label{fig:Test3_convergence}
\end{figure}

\subsubsection{Tests taken form \cite{Vicsek_Riemann}}

\paragraph{Test 4: 1-shock + 2-shock}

This test consists on a 1-shock with speed $\sigma_1 = -0.24252576951298943$ followed by a 2-shock with speed $\sigma_2 =  0.3627332052302428$ linked through an intermediate state $U_{m} = (\rho_m, \theta_m)^{T} = (2.580231243693331, 1.5098899053692902)^{T}$. In Figure \ref{fig:Test4} the Godunov schemes without DR propagate with wrong speeds and the intermediate density is clearly overestimated. With DR, both schemes capture the two shocks sharply and at the correct positions, with only very small oscillations in the intermediate state. The HLL viscous scheme also provides correct shock positions and a smooth intermediate state, but the shocks are noticeably more diffused than with the DR schemes.

\begin{figure}[htpb]
	\begin{subfigure}[h]{0.5\textwidth}
		\centering
			\includegraphics[width=1.1\linewidth]{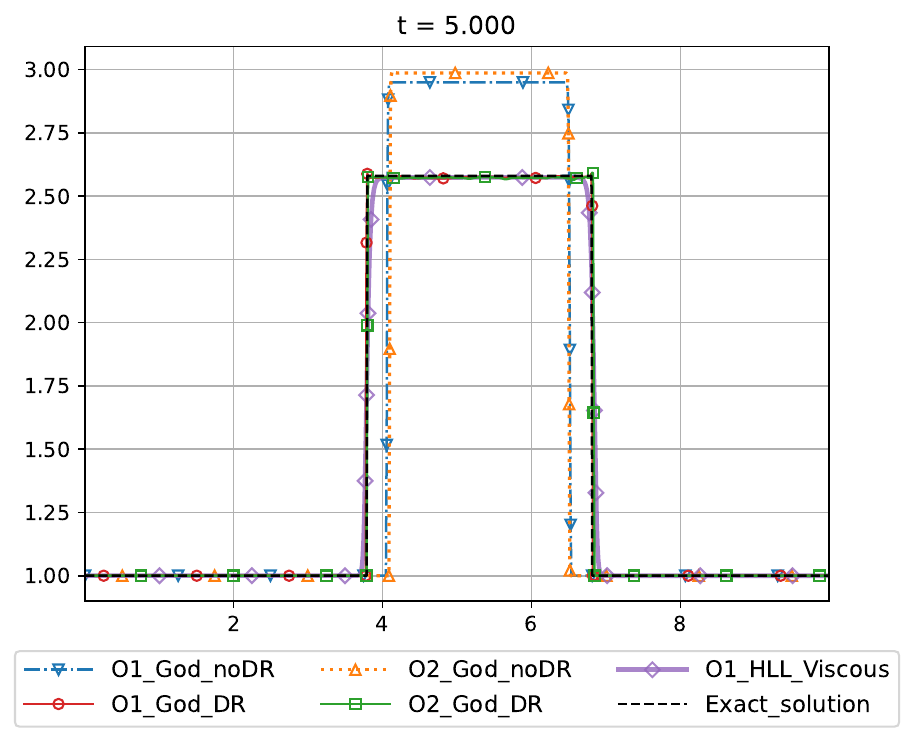}
            \caption{Variable $\rho$}
	\end{subfigure}
	\begin{subfigure}[h]{0.5\textwidth}
		\centering
			\includegraphics[width=1.1\linewidth]{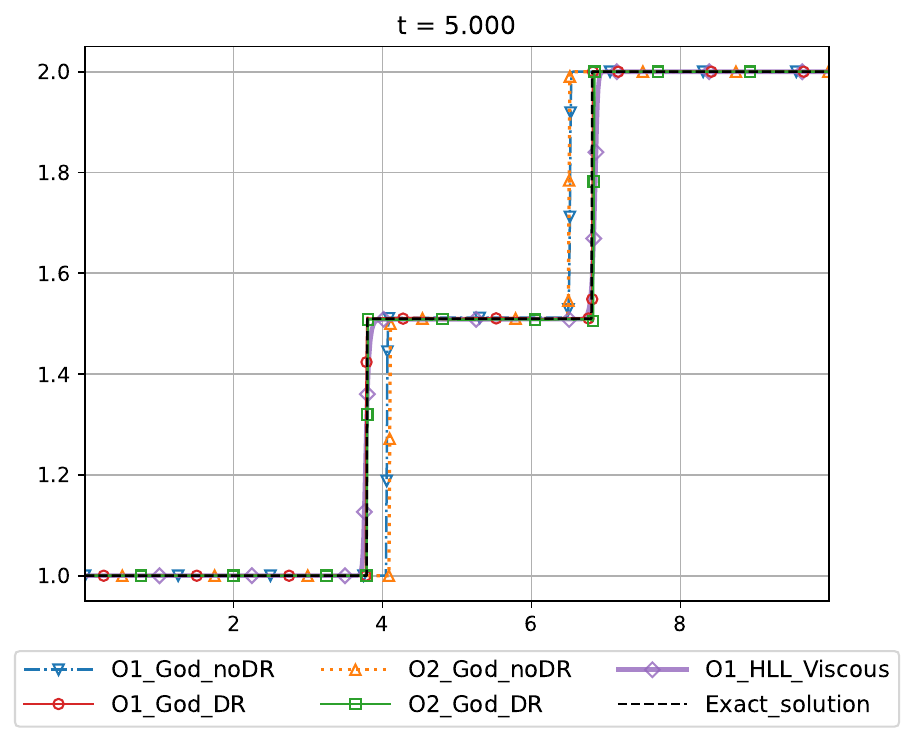}
            \caption{Variable $\theta$}
	\end{subfigure}
	\caption{Numerical results for test 4: 1-shock + 2-shock.}
	\label{fig:Test4}
\end{figure}

\paragraph{Test 5: 1-rarefaction + 2-rarefaction}

This test consists on a 1-rarefaction followed by a 2-rarefaction linked through an intermediate state $U_{m} = (\rho_m, \theta_m)^{T} = (0.21862670530690873, 1.335917425836824)^{T}$. In Figure \ref{fig:Test4}, since the solution contains no shocks, the in-cell discontinuous reconstruction is never activated, so the schemes with and without DR give identical results. The second-order schemes reproduce the edges of the rarefactions more accurately than the first-order ones. The HLL viscous scheme behaves similarly to the first-order Godunov scheme, although it is slightly more diffusive. It also shows a small deviation in the intermediate state.

\begin{figure}[htpb]
	\begin{subfigure}[h]{0.5\textwidth}
		\centering
			\includegraphics[width=1.1\linewidth]{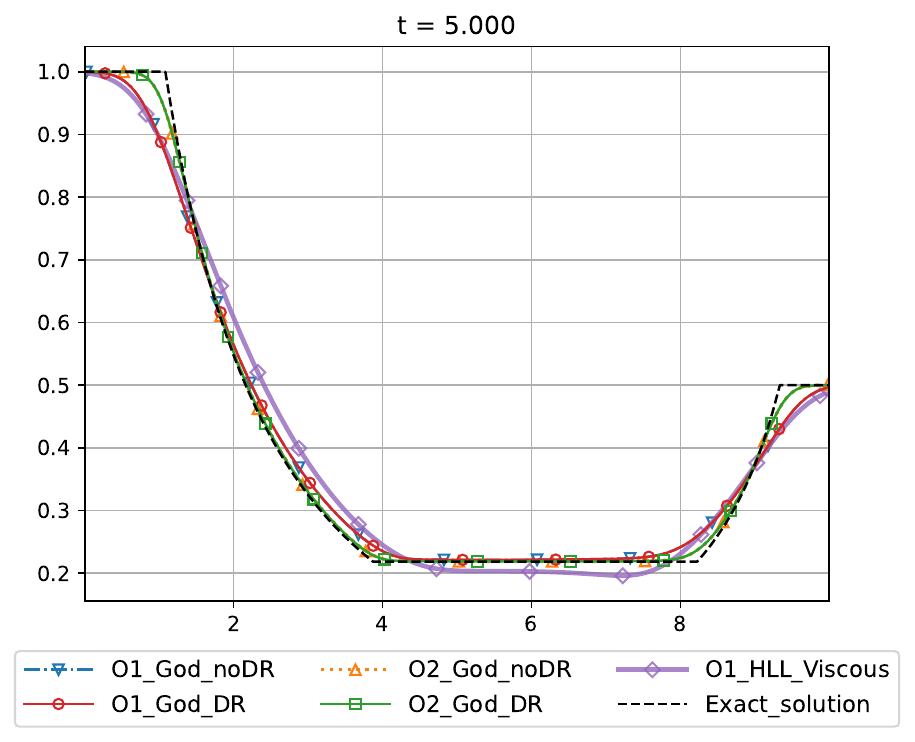}
            \caption{Variable $\rho$}
	\end{subfigure}
	\begin{subfigure}[h]{0.5\textwidth}
		\centering
			\includegraphics[width=1.1\linewidth]{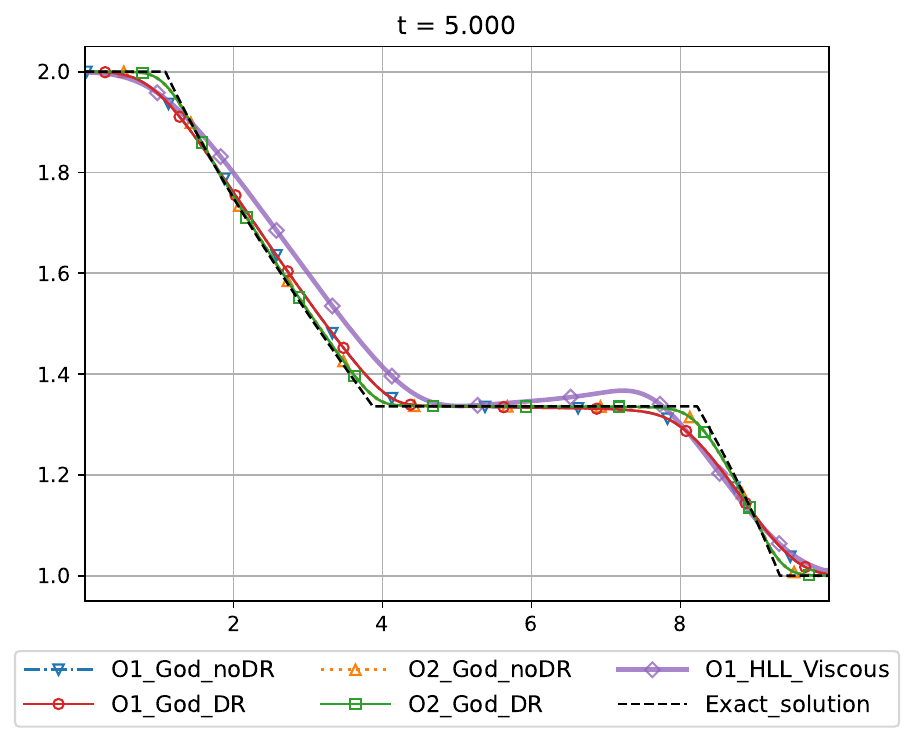}
            \caption{Variable $\theta$}
	\end{subfigure}
	\caption{Numerical results for test 5: 1-rarefaction + 2-rarefaction.}
	\label{fig:Test5}
\end{figure}

\paragraph{Test 6: 1-shock + 2-rarefaction}

This test consists on a 1-shock with speed $\sigma_1 = -0.8665350302618553$ followed by a 2-rarefaction linked through an intermediate state $U_{m} = (\rho_m, \theta_m)^{T} = (2.64567223988045, 2.3314594288198824)^{T}$. In Figure \ref{fig:Test6} the Godunov schemes without DR compute a shock with the wrong speed and an incorrect intermediate state. With DR, the shock and the intermediate state are correctly captured, and the rarefaction is well resolved, especially at second-order. The HLL viscous scheme recovers the correct shock position, but the discontinuity is smeared, the intermediate state is slightly underestimated and the rarefaction is more diffusive.

\begin{figure}[htpb]
	\begin{subfigure}[h]{0.5\textwidth}
		\centering
			\includegraphics[width=1.1\linewidth]{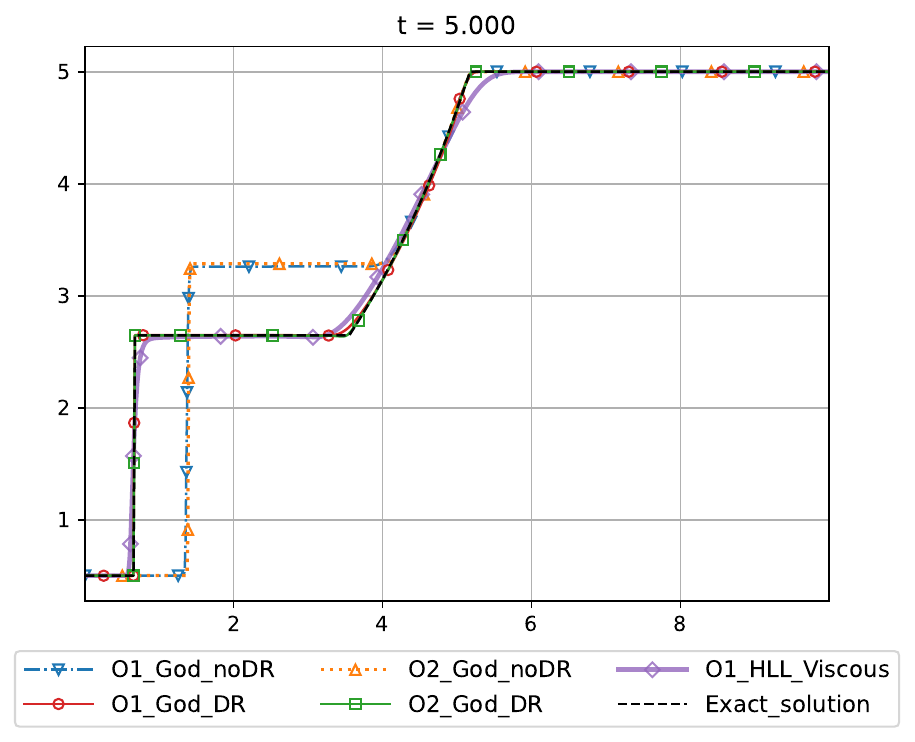}
            \caption{Variable $\rho$}
	\end{subfigure}
	\begin{subfigure}[h]{0.5\textwidth}
		\centering
			\includegraphics[width=1.1\linewidth]{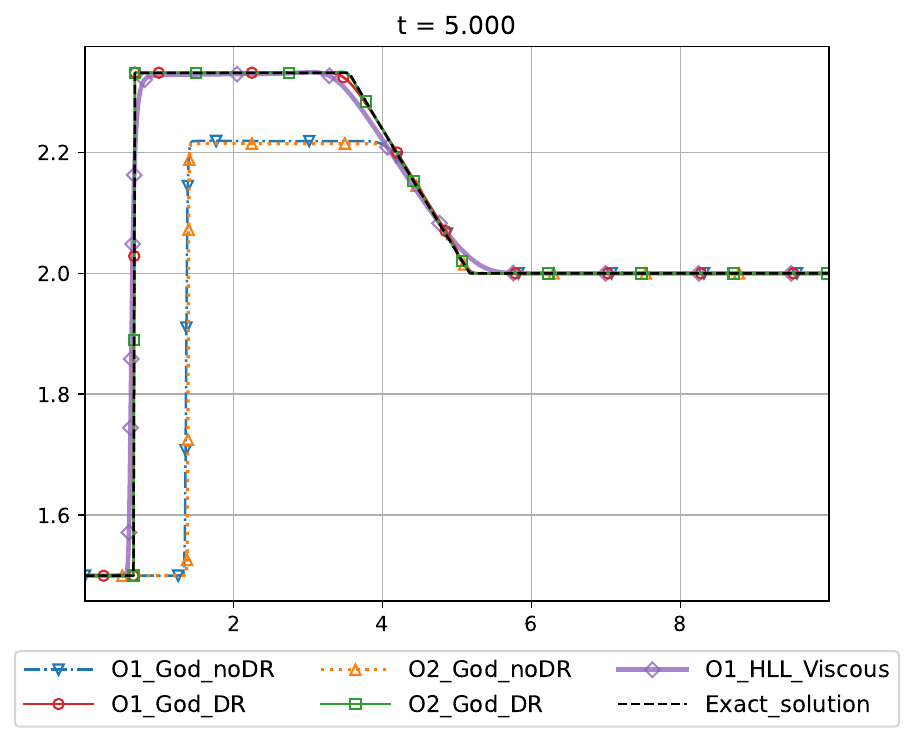}
            \caption{Variable $\theta$}
	\end{subfigure}
	\caption{Numerical results for test 6: 1-shock + 2-rarefaction.}
	\label{fig:Test6}
\end{figure}

\section{Conclusions}

In this work, we have developed a Godunov scheme for the SOH model (\ref{pde}), based on the resolution of the Riemann problem given in \cite{Vicsek_Riemann}. Since the system is nonconservative, the standard Godunov scheme does not converge to the correct shock wave solutions: the numerical shocks travel with wrong speeds and the intermediate states are not correctly recovered, even for a stationary shock. To overcome this difficulty, we have introduced an in-cell discontinuous reconstruction operator into the Godunov scheme. Although this technique has already been used in previous works \cite{chalons2020path, pimentelgarcia2022cell, pimentel2024drmood}, it had always been applied within the framework of path-conservative schemes. To the best of our knowledge, this is the first time it is used outside this framework, here combined with a Godunov scheme based on the exact solution of the Riemann problem. The resulting methods have been implemented at first and second order in space and time. The numerical results, compared with the exact solutions and with the viscous Godunov-type scheme of \cite{compain:hal-05490904}, show the relevance of the in-cell discontinuous reconstruction. With it, the shocks are located correctly, the intermediate states are well recovered, and the shocks are captured sharply, with no numerical viscosity near the discontinuities. This is a clear improvement with respect to the viscous Godunov-type scheme, which smears the shocks over several cells. Finally, a mesh refinement study shows numerically that the schemes with in-cell discontinuous reconstruction converge to the exact solution.

%%%%%%%%%%%%%%%%%%%%%%%%%%%%%%%%%%%%%%%%%%%%%%%%%%%%%%%%%%%
%%%%%%%%%%%%%%%%%%%%%%%%%%%%%%%%%%%%%%%%%%%%%%%%%%%%%%%%%%%
%%%%%%%%%%%%%%%%%%%%%%%%%%%%%%%%%%%%%%%%%%%%%%%%%%%%%%%%%%%

\appendix

\section{Integrals}

\subsection{Primitive function for the jump conditions}\label{appendix_1}
In this appendix, we derive an explicit formulation for the primitive function needed in the evaluation of the Rankine-Hugoniot relations \eqref{G_sigma} for shock wave solutions.
\begin{eqnarray*}
G_\sigma(\theta) & = & \int \frac{- \sigma + c \cos(\theta) - \frac{\lambda \sin(\theta)^2}{\cos(\theta) - \sigma}}{(\cos(\theta) - \sigma)^2}\, d\theta\\
& = & - \frac{\lambda \sin(\theta)}{2(\cos(\theta) - \sigma)^2} + 
\frac{\sigma (2c + \lambda - 2) \sin(\theta)}{2 (1 - \sigma^2) (\cos(\theta) - \sigma)}
+  \frac{2c + \lambda - 2\sigma^2}{2 (1 - \sigma^2)} J(\theta, \sigma)
\end{eqnarray*}
with
$$
J(\theta, \sigma) = \begin{cases}
\frac{2}{\sqrt{1 - \sigma^2}}\mathrm{arctanh}\left( \sqrt{\frac{1 + \sigma}{1 - \sigma}}t \right)
  & \text{if $|\sigma| < 1$;}\\
  -\frac{2}{\sqrt{\sigma^2 -1}}\arctan\left( \sqrt{\frac{\sigma + 1}{\sigma-1}}t\right)
  & \text{if $\sigma > 1$;}\\
    \frac{2}{\sqrt{\sigma^2 -1}}\arctan\left( \sqrt{\frac{-\sigma - 1}{1 - \sigma}}t\right)
  & \text{if $\sigma < -1$;}
  \end{cases}
  $$
  with 
  $$
  t = \tan\left( \frac{\theta}{2}\right).$$

\subsection{Primitive function for the integral curve}\label{appendix_2}
Now, the purpose of this appendix is to give the expression of the two functions whose intersection gives the intermediate value $\theta_m$:
$$
\mathcal{F}_-(\theta) = \begin{cases}
\rho_L\exp(F_-(\theta)-F_-(\theta_L)) & \text{if $\theta \leq \theta_L$;}\\
\displaystyle \rho_L \left(\frac{\sigma_- -\cos(\theta_L)}{\sigma_- -\cos(\theta_R)} \right) & \text{otherwise;}
\end{cases}
$$
and 
$$
\mathcal{F}_+^b(\theta) = \begin{cases}
\displaystyle\rho_R \left(\frac{\sigma_+ -\cos(\theta)}{\sigma_+ -\cos(\theta_L)} \right) & \text{if $\theta \leq \theta_R$;}\\
\rho_R\exp(F_+(\theta)-F_+(\theta_R)) & \text{otherwise.}
\end{cases}
$$
where     \begin{equation}
        F_\pm(x) = -\frac{a}{\lambda} \log \left| \sqrt{a^2\cos^2(x) + \lambda\sin^2(x)} \mp a\cos(x) \right| \pm \frac{\sqrt{\lambda - a^2}}{\lambda} \arcsin\left( \cos(x) \sqrt{\frac{\lambda - a^2}{\lambda}} \right) 
        \end{equation}
        and
        $$
a^2 = \left(\frac{1-c}{2}\right)^2.
$$
  
\section*{Acknowledgments}

Christophe Berthon and Marie Compain acknowledge the support of the ANR project COOKIE ANR-25-CE40-5565.
Marie Compain also acknowledges the support of the CNRS IRP PICASSO.
Ernesto Pimentel-García and Carlos Parés have been supported by projects PID2022-137637NB-C21, PID2025-170651NB-C21 and PCI2024-155061-2 funded by \\MCIN/AEI/10.13039/501100011033/ and ERDF A way of making Europe. Both have also been supported by the Grant PPRO-FQM216-2023 funded by Consejería de Universidad, Investigación e Innovación and by ERDF Andalusia Program 2021-2027.

\bibliographystyle{plain}
\bibliography{references}

\end{document}